\documentclass[12pt]{article}
\usepackage{amsmath,amssymb,amscd}
\title{Gorenstein Biliaison and ACM Sheaves}
\author{Marta Casanellas\footnote{Research of the first author partially
supported by the Secretaria de Estado de Educacion y Universidada of Spain and
the European Social Funding.}
\\ Robin Hartshorne}
\date{}

\begin{document}

\maketitle

\begin{abstract}
Let $X$ be a normal arithmetically Gorenstein scheme in ${\mathbb P}^n$.  We give
a criterion for all codimension two ACM subschemes of $X$ to be in the same
Gorenstein biliaison class on $X$, in terms of the category of ACM sheaves on
$X$.  These are sheaves that correspond to the graded maximal Cohen--Macaulay
modules on the homogeneous coordinate ring of $X$.  Using known results on MCM 
modules, we are able to determine the Gorenstein biliaison classes of codimension
two subschemes of certain varieties, including the nonsingular quadric surface in
${\mathbb P}^3$, and the cone over it in ${\mathbb P}^4$.  As an application we
obtain a new proof of some theorems of Lesperance about curves in ${\mathbb
P}^4$, and answer some questions be raised.
\end{abstract}

\setcounter{section}{-1}
\section{Introduction}
\label{sec0}

Liaison has become an established technique in algebraic geometry.  See, for
example, the excellent book of Migliore \cite{M} for an introduction and more than
200 references.

More recently many people have studied Gorenstein liaison and the important open
question, whether every arithmetically Cohen--Macaulay (ACM) subscheme of
${\mathbb P}^n$ is in the Gorenstein liaison class of a complete intersection
(see, for example, \cite{KMMNP}).

In this paper we study Gorenstein biliaison, defined as the equivalence relation
generated by elementary Gorenstein biliaisons (see \cite[5.4.7]{M},
\cite[1.1]{SEG}, \cite[\S 3]{GDB}) on a normal arithmetically Cohen--Macaulay
(ACM) projective scheme $X$.  We give a new approach to biliaison by relating it
to the category of ACM sheaves on $X$.  Note in this paper biliaison is not a
synonym for even liaison.

Our main theorems $(4.2)$, $(4.3)$ give a criterion for all codimension two ACM
subschemes of $X$ to be in the same Gorenstein biliaison class, in terms of two
conditions on the category of arithmetically Cohen--Macaulay (ACM) sheaves on
$X$.  These sheaves correspond to maximal Cohen--Macaulay (MCM) modules on the
homogeneous coordinate ring $S$ of $X$.  Theorem $(4.7)$ gives a criterion for
the biliaison class of curves on a three-fold to be determined by their Rao
modules.  Then making use of previous work on MCM modules over local rings (see
for example the book of Yoshino \cite{Y}) we are able to conclude results about
Gorenstein biliaison of points and curves on certain projective surfaces and
three-folds.  In particular, we show that all zero-schemes on a non-singular
quadric surface in ${\mathbb P}^3$ $(5.1)$ or on a cubic scroll in ${\mathbb P}^4$
$(5.4)$ are in the same Gorenstein biliaison class.  On the cone in ${\mathbb
P}^4$ over a non-singular quadric surface in ${\mathbb P}^3$, we show that two
curves are in the same Gorenstein biliaison class if and only if their Rao modules
are isomorphic $(6.2)$.  This allows us to give a new proof and extend a theorem
of Lesperance concerning biliaison of curves that are disjoint unions of two
plane curves in ${\mathbb P}^4$ $(6.4)$, $(6.5)$.

Our hope is that a better understanding of biliaison of codimension two
subschemes of hypersurfaces will lead to more insight into the problems of
biliaison of codimension three subschemes of projective space.  In a separate
article \cite{CDH} we discuss Gorenstein liaison (as opposed to biliaison) and
give an analogous criterion in terms of ACM sheaves for all codimension two ACM
schemes to be in the same $G$-liaison class.

Sections 1 and 2 are introductory.  Section 3 contains a criterion for
biliaison.  Section 4 contains the main theorems and their proofs.  Sections 5
and 6 contain applications to particular surfaces and three-folds in projective
space.

\bigskip
\noindent
{\bf Acknowledgement:}  The first author would like to thank the Department of
Mathematics at the University of California at Berkeley for the warm welcome
during the period this research was done.  

\section{Biliaison}
\label{sec1}

Let $V_1,V_2$ be equidimensional closed subschemes of dimension $r$, without
embedded components, of
${\mathbb P}_k^n$, the $n$-dimensional projective space over an algebraically
closed field
$k$.

\bigskip
\noindent
{\bf Definition.} 
A closed subscheme $X \subseteq {\mathbb P}_k^n$ is {\em arithmetically
Cohen--Macaulay} (ACM) if its homogeneous coordinate ring $S(X) =
k[x_0,\dots,x_n]/I_X$ (where $I_X$ is the saturated ideal of $X$) is a
Cohen--Macaulay ring.  This is equivalent to saying $H_*^1({\mathcal
I}_{X,{\mathbb P}^n}) = 0$ and $H_*^i({\mathcal O}_X) = 0$ for $0 < i < \dim X$. 
(For any coherent sheaf ${\mathcal F}$ we denote by $H_*^i({\mathcal F})$ the sum
$\oplus_{\ell \in {\mathbb Z}} H^i({\mathcal F}(\ell))$.)

\bigskip
\noindent
{\bf Definition.}  We say $V_2$ is obtained by an {\em elementary biliaison of
height}
$m$ from
$V_1$ if there exists an ACM scheme $Y \subseteq {\mathbb P}^n$, of dimension
$r+1$ containing $V_1$ and $V_2$, such that $V_2 \sim V_1 + mH$ on $Y$.  (Here
$\sim$ means linear equivalence of divisors on $Y$ in the sense of \cite{GD} and
\cite{GDB}.)  The equivalence relation generated by elementary biliaisons is
called simply {\em biliaison}.  If $Y$ is a complete intersection scheme in
${\mathbb P}^n$, we speak of $CI$-{\em biliaison}.  If $Y$ is ACM and satisfies
$G_0$ (Gorenstein in codimension zero) we speak of {\em Gorenstein biliaison} or
$G$-{\em biliaison}.  If $V_1,V_2,Y$ are contained in some projective scheme $X
\subseteq {\mathbb P}^n$, we speak of {\em biliaison} (resp.~$CI$-{\em biliaison},
$G$-{\em biliaison}) {\em on} X.

\bigskip
Note that $CI$-biliaison in ${\mathbb P}^n$ is equivalent to even $CI$-liaison (in
the sense of \cite[5.1.2]{M}) \cite[4.4]{GD}.  However, on a projective scheme
$X$, the notions of $CI$-biliaison, $G$-biliaison, and even $G$-liaison on $X$ (in
the sense of \cite[5.1.2]{M}) are all three distinct.  So we emphasize that the
word ``biliaison'' is not a synonym for ``even liaison''.  Recall, however, that
every $G$-biliaison is an even $G$-liaison \cite[3.6]{GDB}.

\bigskip
\noindent
{\bf Examples 1.1.}  If $X$ is a non-singular quadric surface in ${\mathbb P}^3$,
then $CI$-biliaison of zero-schemes on $X$ preserves parity of the length of the
zero-scheme.  But on a line in $X$ one can make a biliaison from one point to two
points.  So $CI$-biliaison and $G$-biliaison  are not equivalent.

If $X$ is a non-singular quadric three-fold in ${\mathbb P}^4$, then $\mbox{Pic }
X = {\mathbb Z}$, so every surface in $X$ is a complete intersection.  Thus
$CI$-biliaison and $G$-biliaison coincide on $X$.  But if one takes a rational
quartic curve $C$ in $X$, and a line $L$ in $X$ meeting $C$ in two points, then
$C \cup L$ is an arithmetically Gorenstein scheme, so $C$ can be linked to $L$ by
one Gorenstein liaison.  Since $L$ can be linked to another line $L'$ by another
liaison, $C$ and $L'$ are evenly $G$-linked, but are not equivalent for
$G$-biliaison, since in this case $G$-biliaison preserves parity of degree.

\bigskip
In this paper we will study biliaison of codimension two subschemes of a
normal ACM projective scheme $X$.  Note that the condition of normality implies
every divisor on $X$ satisfies $G_0$, so any biliaison of codimension two
subschemes of
$X$ is a $G$-biliaison, and therefore also an even $G$-liaison \cite[3.6]{GDB}.

For a curve $C$ in ${\mathbb P}^n$, we define as usual the {\em Rao module} of
$C$ to be $M(C) = H_*^1({\mathcal I}_{C,{\mathbb P}^n}) = \oplus_{\ell \in
{\mathbb Z}} H^1({\mathcal I}_{C,{\mathbb P}^n}(\ell))$.  If $C$ is contained in
an ACM scheme $X$ of dimension $\ge 2$, note that the Rao module $M(C)$ can also
be computed as $H_*^1({\mathcal I}_{C,X})$.

\bigskip
\noindent
{\bf Proposition 1.2.}  {\em If curves $C,C'$ in ${\mathbb P}^n$ are equivalent
for biliaison, then $M(C') \cong M(C)(h)$ for some integer $h$.}

\bigskip
\noindent
{\em Proof.}  It is enough to check this for one elementary biliaison.  So
suppose $C' \sim C + mH$ on an ACM surface $Y$.  Then we can compute $M(C) =
H_*^1({\mathcal I}_{C,Y})$ and $M(C') = H_*^1({\mathcal I}_{C',Y})$.  But
${\mathcal I}_{C',Y} \cong {\mathcal I}_{C,Y}(-m)$, so the two Rao modules are
isomorphic up to a shift by $-m$.

\section{ACM sheaves}
\label{sec2}

Let $X$ be an ACM subscheme of ${\mathbb P}^N$.

\bigskip
\noindent
{\bf Definition.}  A coherent sheaf ${\mathcal E}$ on $X$ is an ACM {\em sheaf}
if it is locally Cohen--Macaulay on $X$ and $H_*^i({\mathcal E}) = \oplus_{\ell
\in {\mathbb Z}} H^i(X,{\mathcal E}(\ell)) = 0$ for $0 < i < n = \dim X$.

\bigskip
\noindent
{\bf Proposition 2.1.}  {\em There is a one-to-one correspondence between {\em
ACM} sheaves on $X$ and graded {\em MCM} (maximal Cohen--Macaulay) modules on
$S(X)$ given by ${\mathcal E} \mapsto H_*^0({\mathcal E}) = E$ and $E \mapsto
E^{\sim}$.}

\bigskip
\noindent
{\em Proof.}  Let ${\mathcal E}$ be an ACM sheaf.  Then $E$ is a finitely
generated graded $S(X)$-module with $H_{\mathfrak m}^0(E) = H_{\mathfrak m}^1(E) =
0$ by construction, where ${\mathfrak m}$ is the irrelevant ideal of $S(X)$. 
Furthermore there are isomorphisms
$H_*^i({\mathcal E}) \cong H_{\mathfrak m}^{i+1}(E)$ for $i \ge 1$.  Thus $E$ is
an MCM module by the local cohomology criterion of depth.

Conversely, if $E$ is graded MCM module, let ${\mathcal E} = E^{\sim}$ be the
associated sheaf.  Since $E$ is an MCM module, ${\mathcal E}$ will be a locally
CM sheaf, and the same isomorphisms as above show $H_*^i({\mathcal E}) = 0$ for $0
< i <
\dim X$.

\bigskip
\noindent
{\bf Remark 2.2.}  In the definition of ACM sheaf, we could omit the requirement
${\mathcal E}$ locally CM, if we add a condition $H^0({\mathcal E}(\ell)) = 0$
for $\ell \ll 0$.  Because then $E = H_*^0({\mathcal E})$ will be a finitely
generated module on $S(X)$, the condition $H_*^i({\mathcal E}) = 0$ for $0 < i <
\dim X$ makes it an MCM module, and then ${\mathcal E} = E^{\sim}$ will be
locally CM.

\bigskip
\noindent
{\bf Proposition 2.3.}  {\em Let $\omega$ be a dualizing sheaf on the {\em ACM}
scheme $X$.  For any coherent sheaf ${\mathcal F}$ denote by ${\mathcal
F}^{\omega}$ the sheaf ${\mathcal H}om({\mathcal F},\omega)$.  Then}
\begin{itemize}
\item[a)] {\em The functor ${\mathcal E} \mapsto {\mathcal E}^{\omega}$ is a
contravariant, exact functor on the category of locally Cohen--Macaulay sheaves
on $X$.}
\item[b)] {\em For any such ${\mathcal E}$, there is a natural isomorphism
${\mathcal E}^{\omega\omega} \cong {\mathcal E}$.}
\item[c)] {\em Serre duality gives $H^i({\mathcal E}^{\omega})$ dual to
$H^{n-i}({\mathcal E})$ for any such ${\mathcal E}$, where $n = \dim X$.}
\item[d)] {\em ${\mathcal E}$ is {\em ACM} if and only if ${\mathcal E}^{\omega}$
is {\em ACM}.}
\end{itemize}

\bigskip
\noindent
{\em Proof.}  If ${\mathcal E}$ is locally CM, then by local
duality at the various points of $X$ we find ${\mathcal E}xt^i({\mathcal
E},\omega) = 0$ for all $i > 0$.  Thus the spectral sequence of local and global
Ext degenerates and we find $\mbox{Ext}^i({\mathcal E},\omega) = H^i({\mathcal
E}^{\omega})$ for all $i$.  Serre duality then says $H^i({\mathcal E}^{\omega})$
is dual to $H^{n-i}({\mathcal E})$, where $n = \dim X$.  Furthermore, ${\mathcal
E}$ locally CM implies ${\mathcal E}^{\omega}$ locally CM (repeat the proof of
\cite[1.13,1.14]{GD} with ${\mathcal E}^{\omega}$ in place of ${\mathcal
E}^{\vee}$).  The functor is exact, because if
\[
0 \rightarrow {\mathcal E}' \rightarrow {\mathcal E} \rightarrow {\mathcal E}''
\rightarrow 0
\]
is an exact sequence of locally CM sheaves, then there is an exact sequence
\[
0 \rightarrow {\mathcal E}''{}^{\omega} \rightarrow {\mathcal E}^{\omega}
\rightarrow {\mathcal E}'{}^{\omega} \rightarrow {\mathcal E}xt^1({\mathcal
E}'',\omega) = 0.
\]
It is well-known that ${\mathcal E}^{\omega\omega} \cong {\mathcal E}$
(see, for example, \cite[1.5]{GDB}).  Finally, Serre duality shows that
${\mathcal E}$ is ACM if and only if ${\mathcal E}^{\omega}$ is ACM.

\bigskip
\noindent
{\bf Corollary 2.4.}  {\em If $X$ is arithmetically Gorenstein, the dual of an
{\em ACM} sheaf is again an {\em ACM} sheaf.}

\bigskip
\noindent
{\em Proof.}  Indeed, ${\mathcal E}^{\omega}$ will be a twist of ${\mathcal
E}^{\vee}$.

\bigskip
\noindent
{\bf Remark 2.5.}  If $Y$ is a codimension one ACM subscheme of the ACM scheme
$X$, then the ideal sheaf ${\mathcal I}_Y$ of $Y$ on $X$ is an ACM sheaf on $X$. 
Conversely, if $X$ is normal and ACM, every rank one ACM sheaf on $X$ is
isomorphic to an almost Cartier divisor, and the effective ones correspond to the
codimension one ACM subschemes of $X$ \cite[2.7]{GD}.

\bigskip
\noindent
{\bf Definition.} An ACM sheaf ${\mathcal E}$ on a normal ACM scheme $X$ is {\em
layered} if there exists a filtration
\[
0 = {\mathcal E}_0 \subseteq {\mathcal E}_1 \subseteq \dots \subseteq {\mathcal
E}_r = {\mathcal E}
\]
whose quotients ${\mathcal E}_i/{\mathcal E}_{i-1}$ are rank $1$ ACM sheaves on
$X$ for $i = 1,\dots,r$.

\bigskip
\noindent
{\bf Definition.}  If ${\mathcal E}$ is a torsion-free sheaf of rank $r$ on a
normal scheme $X$ we define its {\em first Chern class} $c_1({\mathcal E})$ to be
the double dual of the highest exterior power $\Lambda^r{\mathcal E}$ of
${\mathcal E}$.  We consider $c_1({\mathcal E})$ as an element of the group
APic $X$ of almost Cartier divisors modulo linear equivalence.  We say
${\mathcal E}$ is {\em orientable} if $c_1({\mathcal E}) \cong {\mathcal
O}_X(\ell)$ is a multiple of the hyperplane class ${\mathcal O}_X(1)$, for some
$\ell \in {\mathbb Z}$.

\bigskip
\noindent
{\bf Definition.}  We say a coherent sheaf ${\mathcal E}$ on $X$ is {\em
dissoci\'e} if ${\mathcal E} \cong \oplus_{i=1}^r {\mathcal O}_X(a_i)$ for some
$a_i \in {\mathbb Z}$.

\section{A criterion for biliaison}
\label{sec3}

We give a criterion for when two codimension $2$ subschemes $V,V'$ of a normal
ACM projective scheme $X$ are in the same $G$-biliaison class.

\bigskip
\noindent
{\bf Theorem 3.1.}  {\em Let $V,V'$ be codimension $2$ subschemes without
embedded components of a normal projective {\em ACM} scheme $X$.  Then $V$ and
$V'$ are in the same $G$-biliaison equivalence class on $X$ if and only if there
exist resolutions on $X$}
\[
0 \rightarrow {\mathcal E} \rightarrow {\mathcal N} \rightarrow {\mathcal
I}_V(a) \rightarrow 0 
\]
\[
0 \rightarrow {\mathcal E}' \rightarrow {\mathcal N} \rightarrow {\mathcal
I}_{V'}(a') \rightarrow 0
\]
{\em with $a,a' \in {\mathbb Z}$, the same coherent sheaf ${\mathcal N}$ in the
middle, and where ${\mathcal E},{\mathcal E}'$ are layered {\em ACM} sheaves
(cf.\ \S 2) on $X$, of the same rank, and the rank $1$ factors of the layerings
of ${\mathcal E}$ and ${\mathcal E}'$ are isomorphic, up to twist, in some order.}

\bigskip
\noindent
{\em Proof.}  First suppose $V$ and $V'$ are equivalent by biliaison.  We proceed
by induction on the number of elementary biliaisons required.

If $V$ and $V'$ are related by a single elementary biliaison, then there is an
ACM divisor $Y$ in $X$, containing $V$ and $V'$, and $V' \sim V + mH$ on $Y$ for
some $m \in {\mathbb Z}$.  This gives an isomorphism of ideal sheaves on $Y$,
${\mathcal I}_{V',Y} \cong {\mathcal I}_{V,Y}(-m)$.  Thus we can write
\begin{eqnarray}
0 \rightarrow {\mathcal I}_Y(-m) \rightarrow {\mathcal I}_V(-m) \rightarrow
&{\mathcal I}_{V,Y}(-m) &\rightarrow 0 \\
&\| \nonumber \\
0 \rightarrow {\mathcal I}_Y \rightarrow {\mathcal I}_{V'} \rightarrow &{\mathcal
I}_{V',Y} &\rightarrow 0.
\end{eqnarray}
Let ${\mathcal F}$ be the fibered sum of ${\mathcal I}_V(-m)$ and ${\mathcal
I}_{V'}$ over ${\mathcal I}_{V',Y}$.  Then we have exact sequences
\begin{equation}
0 \rightarrow {\mathcal I}_Y(-m) \rightarrow {\mathcal F} \rightarrow {\mathcal
I}_{V'} \rightarrow 0
\end{equation}
\begin{equation}
0 \rightarrow {\mathcal I}_Y \rightarrow {\mathcal F} \rightarrow
{\mathcal I}_V(-m) \rightarrow 0.
\end{equation}
Thus the condition of the theorem is satisfied, with ${\mathcal E},{\mathcal E}'$
being the rank $1$ ACM sheaves ${\mathcal I}_Y,{\mathcal I}_Y(-m)$.

Now suppose that $V'$ and $V''$ are related by $r-1$ elementary biliaisons.  By
the induction step, there are sequences
\begin{equation}
0 \rightarrow {\mathcal E}' \rightarrow {\mathcal N} \rightarrow {\mathcal
I}_{V'}(a') \rightarrow 0
\end{equation}
\begin{equation}
0 \rightarrow {\mathcal E}'' \rightarrow {\mathcal N} \rightarrow
{\mathcal I}_{V''}(a'') \rightarrow 0
\end{equation}
with ${\mathcal E}',{\mathcal E}''$ of rank $r-1$, layered, with isomorphic
factors up to twist and order.

Twisting $(3)$ by $a'$ and combining with $(5)$, taking fibered sums as above, we
get new sequences
\begin{equation}
0 \rightarrow {\mathcal I}_Y(-m+a') \rightarrow {\mathcal G} \rightarrow
{\mathcal N} \rightarrow 0
\end{equation}
\begin{equation}
0 \rightarrow {\mathcal E}' \rightarrow {\mathcal G} \rightarrow
{\mathcal F}(a') \rightarrow 0.
\end{equation}
Now composing with the maps ${\mathcal N} \rightarrow {\mathcal I}_{V''}(a'')$
from $(6)$ and ${\mathcal F}(a') \rightarrow {\mathcal I}_V(-m+a')$ from $(4)$ and
using the snake lemma, we get
\begin{equation}
0 \rightarrow {\mathcal H}'' \rightarrow {\mathcal G} \rightarrow {\mathcal
I}_{V''}(a'') \rightarrow 0
\end{equation}
\begin{equation}
0 \rightarrow {\mathcal H} \rightarrow {\mathcal G} \rightarrow {\mathcal
I}_V(-m+a') \rightarrow 0
\end{equation}
where ${\mathcal H}$ and ${\mathcal H}''$ are extensions
\begin{equation}
0 \rightarrow {\mathcal E}' \rightarrow {\mathcal H} \rightarrow {\mathcal
I}_Y(a') \rightarrow 0
\end{equation}
\begin{equation}
0 \rightarrow {\mathcal I}_Y(-m+a') \rightarrow {\mathcal H}'' \rightarrow
{\mathcal E}'' \rightarrow 0.
\end{equation}
Thus ${\mathcal H},{\mathcal H}''$ are layered ACM sheaves with isomorphic
factors, up to twist and order, and so sequences $(9)$, $(10)$ satisfy the
condition of the theorem.

\bigskip
For the reverse implication of the theorem, we will need several lemmas.

\bigskip
\noindent
{\bf Lemma 3.2.}  {\em Let $X$ be an integral projective scheme satisfying the
condition $S_2$ of Serre.  Let ${\mathcal E}$ be a torsion-free coherent sheaf on
$X$, locally free in codimension $1$.  Let $W \subseteq H^0({\mathcal E})$ be a
subspace, and let ${\mathcal E}_0$ be the subsheaf generated by $W$.  Then the
following conditions are equivalent.}

\begin{itemize}
\item[(i)] {\em There is an $s \in W$ such that ${\mathcal E}' = {\mathcal
E}/(s)$ is torsion-free and locally free in codimension $1$.}
\item[(ii)] (a) {\em for all $x \in X$ of codimension $1$,
$\mbox{\em rank}({\mathcal E}_0 \otimes k(x) \stackrel{\sigma_x}{\rightarrow}
{\mathcal E} \otimes k(x)) \ge 1$ and}
\item[] (b) {\em either $\mbox{{\em rank} } {\mathcal E}_0 \ge 2$ or ${\mathcal
E}_0
\cong {\mathcal O}_X$ and ${\mathcal E}/{\mathcal E}_0$ is torsion-free and
locally free in codimension $1$.}
\end{itemize}

\bigskip
\noindent
{\em Proof.}  \cite[2.6]{RTLRP}.

\bigskip
\noindent
{\bf Lemma 3.3 (The case $\mbox{rank } {\mathcal N} = 2$).}  {\em With the
hypotheses of the theorem suppose given}
\begin{equation}
0 \rightarrow {\mathcal L} \rightarrow {\mathcal N} \rightarrow {\mathcal I}_V(a)
\rightarrow 0
\end{equation}
\begin{equation}
0 \rightarrow {\mathcal L}' \rightarrow {\mathcal N} \rightarrow {\mathcal
I}_{V'}(a') \rightarrow 0
\end{equation}
{\em with $\mbox{{\em  rank} } {\mathcal N} = 2$ and ${\mathcal L},{\mathcal L}'$
{\em ACM} sheaves of $\mbox{{\em rank} } 1$, isomorphic up to twist.  Then either
$V = V'$ or $V$ and $V'$ are related by a single elementary biliaison on $X$.}

\bigskip
\noindent
{\em Proof.}  Make a diagram
\[
\begin{array}{ccccccccc}
& & & &0 & &0 \\
& & & &\downarrow & &\downarrow \\
& & & &{\mathcal L}' &= &{\mathcal L}' \\
& & & &\downarrow & &\ \ \,\downarrow \alpha \\
0 &\rightarrow &{\mathcal L} &\rightarrow &{\mathcal N} &\rightarrow &{\mathcal
I}_V(a) &\rightarrow &0 \\
& &\| & &\downarrow & &\downarrow \\
0 &\rightarrow &{\mathcal L} &\rightarrow &{\mathcal I}_{V'}(a') &\rightarrow
&{\mathcal R} &\rightarrow &0 \\
& & & &\downarrow & &\downarrow \\
& & & &0 & &0
\end{array}
\]
We start by considering the composed map $\alpha: {\mathcal L}' \rightarrow
{\mathcal I}_V(a)$.  If $\alpha = 0$, then ${\mathcal L}'$ maps to ${\mathcal L}$
with torsion free cokernel, so ${\mathcal L} \cong {\mathcal L}'$, $a = a'$, $V =
V'$.  Otherwise $\alpha$ must be injective, and we get the diagram shown.  Now
the map $\alpha$ composed with the inclusion ${\mathcal I}_V(a) \subseteq
{\mathcal O}_X(a)$ gives the ideal of an ACM divisor $Y$ on $X$, so that
${\mathcal L}' = {\mathcal I}_Y(a)$. 
Then $V$ is contained in $Y$ and ${\mathcal R} \cong {\mathcal I}_{V,Y}(a)$. 
From the bottom row of the diagram we see also that ${\mathcal L} = {\mathcal
I}_Y(a')$ and 
${\mathcal R}
\cong {\mathcal I}_{V',Y}(a')$.  Thus ${\mathcal I}_{V,Y}(a) \cong {\mathcal
I}_{V',Y'}(a')$, so
$V' \sim V + (a'-a)H$ on $Y$.

\bigskip
\noindent
{\bf Lemma 3.4.}  {\em Let $V$ be a codimension $2$ subscheme of the normal {\em
ACM} scheme $X$ and suppose given a sequence}
\[
0 \rightarrow {\mathcal E} \rightarrow {\mathcal N} \rightarrow {\mathcal I}_V(a)
\rightarrow 0
\]
{\em where ${\mathcal N}$ is coherent and ${\mathcal E}$ is a layered {\em ACM}
sheaf with $\mbox{{\em rank} } 1$ {\em ACM} factors ${\mathcal L}_i$, $i =
1,\dots,r$.  Then for any choice of $b_i \gg 0$ and $s_i \in
\mbox{\em Hom}({\mathcal L}_i,{\mathcal N}(b_i))$ sufficiently general, we will
get a sequence}
\[
0 \rightarrow \oplus {\mathcal L}_i(-b_i) \stackrel{\oplus s_i}{\longrightarrow}
{\mathcal N} \rightarrow {\mathcal I}_{V'}(a') \rightarrow 0
\]
{\em for some codimension $2$ subscheme $V'$ of $X$ in the same $G$-biliaison
class as $V$.}

\bigskip
\noindent
{\em Proof.}  We will split off the factors of ${\mathcal E}$ one by one.  Since
${\mathcal E}$ is layered, we can write
\[
0 \rightarrow {\mathcal E}' \rightarrow {\mathcal E} \rightarrow {\mathcal L}
\rightarrow 0
\]
where ${\mathcal E}'$ is layered of rank $r-1$, and ${\mathcal L}$ is one of the
factors.  We write a diagram
\[
\begin{array}{ccccccccc}
& &0 & &0 \\
& &\downarrow & &\downarrow \\
& &{\mathcal E}' &= &{\mathcal E}' \\
& &\downarrow & &\downarrow \\
0 &\rightarrow &{\mathcal E} &\rightarrow &{\mathcal N} &\rightarrow &{\mathcal
I}_V(a) &\rightarrow &0 \\
& &\downarrow & &\downarrow & &\| \\
0 &\rightarrow &{\mathcal L} &\rightarrow &{\mathcal F} &\rightarrow &{\mathcal
I}_V(a) &\rightarrow &0 \\
& &\downarrow & &\downarrow \\
& &0 & &0
\end{array}
\]
where ${\mathcal F} = {\mathcal N}/{\mathcal E}'$ is of rank $2$.  Choose $b \gg
0$ large enough so that the sheaf ${\mathcal H}om({\mathcal L},{\mathcal N}(b))$
is generated by global sections, and let $W = H^0({\mathcal H}om({\mathcal
L},{\mathcal F}(b)))$ be the image of $H^0({\mathcal H}om({\mathcal L},{\mathcal
N}(b)))$.  There is an exact sequence
\[
0 \rightarrow {\mathcal H}om({\mathcal L},{\mathcal E}') \rightarrow {\mathcal
H}om({\mathcal L},{\mathcal N}) \rightarrow {\mathcal H}om({\mathcal L},{\mathcal
F}) \rightarrow {\mathcal E}xt^1({\mathcal L},{\mathcal E}') \rightarrow \dots\,.
\]
Since $X$ is normal, ${\mathcal L}$ is locally free in codimension $1$, so the
sheaf ${\mathcal E}xt^1({\mathcal L},{\mathcal E}')$ has support in codimension
$\ge 2$.  It follows that $W$ generates a subsheaf ${\mathcal F}_0$ of ${\mathcal
H}om({\mathcal L},{\mathcal F}(b))$ that is equal to ${\mathcal H}om({\mathcal
L},{\mathcal F}(b))$ in codimension $1$ and therefore has rank $= 2$.  Thus we
can apply $(3.2)$ and find that for any sufficiently general $s \in
\mbox{Hom}({\mathcal L},{\mathcal N}(b))$, the composed map ${\mathcal L}(-b)
\rightarrow {\mathcal F}$ will have cokernel torsion-free and locally free in
codimension $1$.  Its first Chern class will be $a+b$, so we get a sequence
\[
0 \rightarrow {\mathcal L}(-b) \stackrel{s}{\rightarrow} {\mathcal F} \rightarrow
{\mathcal I}_{V'}(a+b) \rightarrow 0
\]
for some codimension $2$ subscheme $V'$.  According to $(3.3)$, $V$ and $V'$ are
related by a single biliaison.

Let ${\mathcal N} \rightarrow {\mathcal I}_{V'}(a+b)$ be the composed map, and
let ${\mathcal R}$ be the kernel.  Then we have a new diagram
\[
\begin{array}{ccccccccc}
& &0 & &0 \\
& &\downarrow & &\downarrow \\
& &{\mathcal E}' &= &{\mathcal E}' \\
& &\downarrow & &\downarrow \\
0 &\rightarrow &{\mathcal R} &\rightarrow &{\mathcal N} &\rightarrow &{\mathcal
I}_{V'}(a+b) &\rightarrow &0 \\
& &\downarrow &{}^s\!\!\!\nearrow &\downarrow & &\| \\
0 &\rightarrow &{\mathcal L}(-b) &\stackrel{s}{\rightarrow} &{\mathcal F}
&\rightarrow &{\mathcal I}_{V'}(a+b) &\rightarrow &0 \\
& &\downarrow & &\downarrow \\
& &0 & &0
\end{array}
\]
and since the map $s$ lifts to ${\mathcal N}$ by construction, and its image in
${\mathcal I}_{V'}(a+b)$ is zero, the first column sequence splits, and we find
${\mathcal R} \cong {\mathcal E}' \oplus {\mathcal L}(-b)$.  Furthermore, the map
${\mathcal L}(-b) \rightarrow {\mathcal N}$ is the chosen one $s$.

We repeat this procedure with a factor of ${\mathcal E}'$.  Continuing in the
same manner, after $r$ steps we obtain the result of the statement, each time
replacing $V$ by something in the same biliaison class.

\bigskip
\noindent
{\em Proof of $3.1$, continued.}  Suppose given two subschemes $V,V'$ and
sequences
\[
0 \rightarrow {\mathcal E} \rightarrow {\mathcal N} \rightarrow {\mathcal I}_V(a)
\rightarrow 0
\]
\[
0 \rightarrow {\mathcal E}' \rightarrow {\mathcal N} \rightarrow {\mathcal
I}_{V'}(a') \rightarrow 0
\]
as in the theorem.  Here ${\mathcal E}$ and ${\mathcal E}'$ are layered ACM
sheaves with factors ${\mathcal L}_1,\dots,{\mathcal L}_r$ isomorphic up to twist
and order.  Applying $(3.4)$ to both of them, we obtain
\[
0 \rightarrow \oplus {\mathcal L}_i(-b_i) \stackrel{\oplus s_i}{\longrightarrow}
{\mathcal N} \rightarrow {\mathcal I}_{V_1}(a_1) \rightarrow 0
\]
\[
0 \rightarrow \oplus {\mathcal L}'_i(-b'_i) \stackrel{\oplus
s'_i}{\longrightarrow} {\mathcal N} \rightarrow {\mathcal I}_{V'_1}(a'_1)
\rightarrow 0
\]
where $V$ and $V_1$ are biliaison equivalent, and $V',V'_1$ are biliaison
equivalent.

But now, since the factors of ${\mathcal E}$ and ${\mathcal E}'$ were isomorphic
up to shift and order, we may assume ${\mathcal L}_i = {\mathcal L}'_i$ for each
$i$.  Also, we can take the $b_i$ in $(3.4)$ sufficiently large arbitrarily, so
we may assume $b_i = b'_i$ for each $i$.  And when we choose the section $s'_i$,
we may assume that they are equal to the $s_i$.  Then $a_1 = a'_1$, $V_1 = V'_1$,
and we are done.

\section{The main theorem}
\label{sec4}

In this section we will show that the property of ACM codimension $2$ subschemes
of $X$ all being in the same biliaison class is equivalent to two statements about
the category of ACM sheaves on $X$.

Let $X$ be a normal ACM scheme in some projective space, of dimension at least
$2$.  Our first condition is

\begin{quote}
\begin{itemize}
\item[(A)] Any two ACM subschemes $V,V'$ of $X$ of codimension $2$ are equivalent
for Gorenstein biliaison.
\end{itemize}
\end{quote}

Our second condition is

\begin{quote}
\begin{itemize}
\item[(B)] Every orientable (cf.~\S 2) ACM sheaf ${\mathcal E}$ on $X$ has a
presentation
\[
0 \rightarrow {\mathcal F}_2 \rightarrow {\mathcal F}_1 \rightarrow {\mathcal E}
\rightarrow 0
\]
where ${\mathcal F}_1$ and ${\mathcal F}_2$ are layered ACM sheaves.
\end{itemize}
\end{quote}

To state our third condition we need some notation.  Let  ${\mathcal M}$ be the
category of layered ACM sheaves on $X$.  Let $G({\mathcal M})$ be the
Grothendieck group of the category ${\mathcal M}$, generated by the objects of
${\mathcal M}$ and with relations ${\mathcal F} - {\mathcal F}' - {\mathcal F}''$
whenever there is an exact sequence $0 \rightarrow {\mathcal F}' \rightarrow
{\mathcal F} \rightarrow {\mathcal F}'' \rightarrow 0$ in ${\mathcal M}$.  We
regard $G({\mathcal M})$ as a ${\mathbb Z}[h]$-module, where $h \cdot {\mathcal
F} = {\mathcal F}(1)$ for any ${\mathcal F} \in {\mathcal M}$.  We consider the
quotient group $G' = G({\mathcal M})/(1-h)G({\mathcal M})$.  Roughly speaking,
$G'$ is the group generated by the rank $1$ ACM sheaves on $X$, with relations
coming from exact sequences in ${\mathcal M}$, and where each sheaf is identified
with all of its twists.  There is a natural homomorphism $c_1: G({\mathcal M})
\rightarrow \mbox{APic } X$ obtained by taking the first Chern class.  This
passes to the quotient to give a map $c_1: G' \rightarrow \mbox{APic }
X/{\mathbb Z} \cdot {\mathcal O}_X(1)$.  Our third condition is

\begin{quote}
\begin{itemize}
\item[(C)] If ${\mathcal E}$ is an orientable layered ACM sheaf on $X$, then its
class in $G'$ is equal to $r \cdot {\mathcal O}$, where $r = \mbox{rank }
{\mathcal E}$.
\end{itemize}
\end{quote}

Note that (C) is equivalent to saying the kernel of the map $c_1: G' \rightarrow
\mbox{APic } X/{\mathbb Z}\cdot{\mathcal O}_X(1)$ is just the subgroup
${\mathbb Z}\cdot {\mathcal O}$ of $G'$.

Our goal is to prove (A) $\Leftrightarrow$ (B) $+$ (C).  We split the proof in
several parts.

\bigskip
\noindent
{\bf Lemma 4.1.}  {\em Let $X$ be a normal {\em ACM} scheme, and let $V,V'$ be
two codimension $2$ subschemes in the same biliaison class.  Suppose also given
exact sequences}
\setcounter{equation}{0}
\begin{equation}
0 \rightarrow {\mathcal E} \rightarrow {\mathcal G} \rightarrow {\mathcal I}_V
\rightarrow 0
\end{equation}
\begin{equation}
0 \rightarrow {\mathcal E}' \rightarrow {\mathcal G}' \rightarrow {\mathcal
I}_{V'}
\rightarrow 0
\end{equation}
{\em with ${\mathcal E},{\mathcal E}',{\mathcal G},{\mathcal G}'$ coherent
sheaves.  Then there are other coherent sheaves ${\mathcal H},{\mathcal
R},{\mathcal S}$, and layered {\em ACM} sheaves ${\mathcal F},{\mathcal F}'$ with
the same {\em rank} $1$ factors up to order and twist, and exact sequences with
suitable twists $a,a'$,}
\begin{equation}
0 \rightarrow {\mathcal R} \rightarrow {\mathcal H} \rightarrow {\mathcal G}(a)
\rightarrow 0
\end{equation}
\begin{equation}
0 \rightarrow {\mathcal E}'(a') \rightarrow {\mathcal R} \rightarrow {\mathcal F}
\rightarrow 0
\end{equation}
\begin{equation}
0 \rightarrow {\mathcal S} \rightarrow {\mathcal H} \rightarrow {\mathcal G}'(a')
\rightarrow 0
\end{equation}
\begin{equation}
0 \rightarrow {\mathcal E}(a) \rightarrow {\mathcal S} \rightarrow {\mathcal F}'
\rightarrow 0.
\end{equation}

\bigskip
\noindent
{\em Proof.}  By $(3.1)$ there are exact sequences
\begin{equation}
0 \rightarrow {\mathcal F} \rightarrow {\mathcal N} \rightarrow {\mathcal I}_V(a)
\rightarrow 0
\end{equation}
\begin{equation}
0 \rightarrow {\mathcal F}' \rightarrow {\mathcal N} \rightarrow {\mathcal
I}_{V'}(a') \rightarrow 0
\end{equation}
with the same coherent sheaf ${\mathcal N}$ in the middle, and layered ACM
sheaves ${\mathcal F},{\mathcal F}'$ with the same rank $1$ factors, up to
twist.  We make fibered sum constructions, as in the proof of $(3.1)$.

From $(1)$ and $(7)$ we obtain
\begin{equation}
0 \rightarrow {\mathcal E}(a) \rightarrow {\mathcal A} \rightarrow {\mathcal N}
\rightarrow 0
\end{equation}
\begin{equation}
0 \rightarrow {\mathcal F} \rightarrow {\mathcal A} \rightarrow {\mathcal G}(a)
\rightarrow 0.
\end{equation}
From $(2)$ and $(8)$ we obtain
\begin{equation}
0 \rightarrow {\mathcal E}'(a') \rightarrow {\mathcal B} \rightarrow {\mathcal N}
\rightarrow 0
\end{equation}
\begin{equation}
0 \rightarrow {\mathcal F}' \rightarrow {\mathcal B} \rightarrow {\mathcal
G}'(a') \rightarrow 0.
\end{equation}
Then from $(9)$ and $(11)$ we obtain
\begin{equation}
0 \rightarrow {\mathcal E}(a) \rightarrow {\mathcal H} \rightarrow {\mathcal B}
\rightarrow 0
\end{equation}
\begin{equation}
0 \rightarrow {\mathcal E}'(a') \rightarrow {\mathcal H} \rightarrow {\mathcal A}
\rightarrow 0.
\end{equation}
Now compose the maps ${\mathcal H} \rightarrow {\mathcal B} \rightarrow {\mathcal
G}'(a')$ and let the kernel be ${\mathcal S}$, to get $(5)$.  Then because of
$(12)$ and $(13)$ we get $(6)$.  Similarly, comparing ${\mathcal H} \rightarrow
{\mathcal A} \rightarrow {\mathcal G}(a)$ we get $(3)$ and $(4)$.

\bigskip
\noindent
{\bf Proposition 4.2.}  {\em Let $X$ be a normal projective {\em ACM} scheme. 
Then {\em (A)} $\Rightarrow$ {\em (B)} $+$ {\em (C)}.}

\bigskip
\noindent
{\em Proof.}  To prove (B), let ${\mathcal E}$ be any orientable ACM sheaf.  Then
we can find a dissoci\'e sheaf ${\mathcal L}$ and an exact sequence
\[
0 \rightarrow {\mathcal L} \rightarrow {\mathcal E} \rightarrow {\mathcal I}_V(b)
\rightarrow 0
\]
for some codimension $2$ subscheme of $V$ \cite[1.12]{RTLRP}.  Furthermore, since
${\mathcal E}$ is an ACM sheaf, $V$ will be an ACM subscheme.  By hypothesis (A)
there is a complete intersection $Z$ of $X$ with two hypersurfaces in the same
biliaison class as $V$.  Let
\[
0 \rightarrow {\mathcal P}_2 \rightarrow {\mathcal P}_1 \rightarrow {\mathcal
I}_Z \rightarrow 0
\]
be its resolution, where ${\mathcal P}_i$ are dissoci\'e sheaves.

Now by the lemma, there are sheaves ${\mathcal H}$, ${\mathcal R}$, ${\mathcal
S}$, ${\mathcal F}$, ${\mathcal F}'$ and exact sequences $(3)$--$(6)$ in changed
notation.  Sequence $(3)$ gives a resolution
\[
0 \rightarrow {\mathcal R} \rightarrow {\mathcal H} \rightarrow {\mathcal E}(a-b)
\rightarrow 0.
\]
Then sequence $(4)$ shows that ${\mathcal R}$ is layered, and sequences $(5)$,
$(6)$ show that ${\mathcal H}$ is layered.  This proves (B).

To prove (C), let ${\mathcal E}$ be an orientable layered ACM sheaf.  Let
${\mathcal L},V,Z,{\mathcal P}_1,{\mathcal P}_2$ be as in the first part of the
proof.  Then the sequences $(3)$--$(6)$ of the lemma show that in the quotient
Grothendieck group $G'$, we have
\[
{\mathcal H} = {\mathcal E} + {\mathcal P}_2 + {\mathcal F} = {\mathcal P}_1 +
{\mathcal L} + {\mathcal F}'.
\]
Since ${\mathcal P}_1$, ${\mathcal P}_2$ and ${\mathcal L}$ are dissoci\'e, and
${\mathcal F}$ and ${\mathcal F}'$ have the same rank $1$ factors, up to twist,
we find ${\mathcal E} = r \cdot {\mathcal O}$ in $G'$, where $r = \mbox{rank }
{\mathcal E}$.

\bigskip
\noindent
{\bf Theorem 4.3.}  {\em Let $X$ be a normal {\em ACM} projective scheme, and
assume furthermore either {\em a)} $X$ is arithmetically Gorenstein, or {\em b)}
$\mbox{{\em APic} } X$ is generated as a monoid by the {\em rank} $1$ {\em ACM}
sheaves.  Then {\em (B) $+$ (C) $\Rightarrow$ (A)}.}

\bigskip
We will need several lemmas before proving the theorem.

\bigskip
\noindent
{\bf Lemma 4.4.}  {\em Under either hypothesis {\em a)} or {\em b)} of $(4.3)$,
the condition {\em (B)} implies}

\begin{quote}
\begin{itemize}
\item[(B${}^*$)] {\em Every orientable {\em ACM} sheaf ${\mathcal E}$ on $X$ has a
resolution}
\[
0 \rightarrow {\mathcal E} \rightarrow {\mathcal G}_1 \rightarrow {\mathcal G}_2
\rightarrow 0
\]
{\em where ${\mathcal G}_1$ and ${\mathcal G}_2$ are layered {\em ACM} sheaves.}
\end{itemize}
\end{quote}

\bigskip
\noindent
{\em Proof.}  a) If $X$ is arithmetically Gorenstein, then the dual of an
orientable ACM sheaf is again an orientable ACM sheaf $(2.4)$.  So we apply (B) to
${\mathcal E}^{\vee}$ and dualize.

b) Suppose that APic $X$ is generated as a monoid by rank $1$ ACM
sheaves.  Then for any ACM sheaf ${\mathcal E}$ we can write $-c_1({\mathcal E})
= \Sigma {\mathcal L}_i$ in APic $X$, where the ${\mathcal L}_i$
are rank
$1$ ACM sheaves.  Then ${\mathcal E}' = {\mathcal E} \oplus \bigoplus {\mathcal
L}_i$ is an orientable ACM sheaf.  Applying (B) to ${\mathcal E}'$ we get
\[
0 \rightarrow {\mathcal F}_2 \rightarrow {\mathcal F}_1 \rightarrow {\mathcal E}'
\rightarrow 0
\]
with the ${\mathcal F}_i$ layered ACM sheaves.  Now composing with the projection
${\mathcal E}' \rightarrow {\mathcal E}$ we get
\[
0 \rightarrow {\mathcal R} \rightarrow {\mathcal F}_1 \rightarrow {\mathcal E}
\rightarrow 0
\]
where ${\mathcal R}$ is an extension
\[
0 \rightarrow {\mathcal F}_2 \rightarrow {\mathcal R} \rightarrow \bigoplus
{\mathcal L}_i \rightarrow 0.
\]
Thus ${\mathcal E}$ has a resolution of the same type, and we have established
the stronger condition

\begin{quote}
\begin{itemize}
\item[(B$'$)] Every ACM sheaf ${\mathcal E}$ on $X$ (not necessarily orientable)
has a resolution of the type (B).
\end{itemize}
\end{quote}

Now we can apply the functor ${\mathcal E} \mapsto {\mathcal E}^{\omega}$ of
$(2.3)$ to the category of all ACM sheaves and obtain the condition

\begin{quote}
\begin{itemize}
\item[(B$'{}^*$)] Every ACM sheaf ${\mathcal E}$ on $X$ has a resolution of type
(B${}^*$).
\end{itemize}
\end{quote}

Since (B$'{}^*$) obviously implies (B${}^*$) we are done.

\bigskip
\noindent
{\bf Lemma 4.5.}  {\em Let $X$ be a normal {\em ACM} projective scheme, and
suppose given a codimension two subscheme $V$ and exact sequences}
\[
0 \rightarrow {\mathcal E}_0 \oplus {\mathcal E}' \rightarrow {\mathcal N}
\rightarrow {\mathcal I}_V(a) \rightarrow 0
\]
\[
0 \rightarrow {\mathcal E}' \rightarrow {\mathcal E} \rightarrow {\mathcal F}
\rightarrow 0
\]
{\em with ${\mathcal E}_0$ and ${\mathcal N}$ coherent, and ${\mathcal
E}',{\mathcal E},{\mathcal F}$ layered {\em ACM} sheaves on $X$.  Then for  any
$b \gg 0$ there exists a sequence}
\[
0 \rightarrow {\mathcal E}_0 \oplus {\mathcal E}'(-b) \rightarrow {\mathcal N}
\rightarrow {\mathcal I}_{V'}(a') \rightarrow 0
\]
{\em with $V'$ in the same biliaison  class as $V$, and with the extra property
that the map ${\mathcal E}'(-b) \rightarrow {\mathcal N}$ in this sequence
extends to a map of ${\mathcal E}(-b) \rightarrow {\mathcal N}$.}

\bigskip
\noindent
{\em Proof.} By induction on the rank of ${\mathcal E}'$.  If ${\mathcal E}'$ has
rank $1$, we can argue as in the proof of $(3.4)$.  Consider the sequence of
sheaves
\[
{\mathcal H}om({\mathcal E},{\mathcal N}) \rightarrow {\mathcal H}om({\mathcal
E}',{\mathcal N}) \rightarrow {\mathcal E}xt^1({\mathcal F},{\mathcal N}).
\]
Since $X$ is normal, ${\mathcal F}$ is locally free in codimension $1$, and so the
${\mathcal E}xt$ sheaf has support in codimension $\ge 2$.  Take any $b \gg 0$ so
that
${\mathcal H}om({\mathcal E},{\mathcal N}(b))$ is generated by global sections,
and let $W
\subseteq H^0({\mathcal H}om({\mathcal E}',{\mathcal N}(b)))$ be the image of
$H^0({\mathcal H}om({\mathcal E},{\mathcal N}(b))$.  Then $W$ satisfies the
conditions of $(3.2)$, and so a general $s \in W$ will give the required map of
${\mathcal E}'(-b) \rightarrow {\mathcal N}$.

Now suppose rank ${\mathcal E}' \ge 2$.  Then we can split off a rank $1$ factor
${\mathcal L}$
\[
0 \rightarrow {\mathcal E}'' \rightarrow {\mathcal E}' \rightarrow {\mathcal L}
\rightarrow 0
\]
and get
\[
0 \rightarrow {\mathcal E}'' \rightarrow {\mathcal E} \rightarrow {\mathcal R}
\rightarrow 0
\]
where
\[
0 \rightarrow {\mathcal L} \rightarrow {\mathcal R} \rightarrow {\mathcal F}
\rightarrow 0.
\]

First apply the splitting technique of $(3.4)$ to get 
\[
0 \rightarrow {\mathcal E}_0 \oplus {\mathcal E}'' \oplus {\mathcal L}(-b_1)
\rightarrow {\mathcal N} \rightarrow {\mathcal I}_{V_1}(a_1) \rightarrow 0
\]
for some $b_1 \gg 0$ and $V_1$ in the same biliaison class as $V$.  Now apply the
induction hypothesis to ${\mathcal E}''$ to find a sequence
\[
0 \rightarrow {\mathcal E}_0 \oplus {\mathcal E}''(-b_2) \oplus {\mathcal
L}(-b_1) \rightarrow {\mathcal N} \rightarrow {\mathcal I}_{V_2}(a_2) \rightarrow
0
\]
for some $b_2 \gg 0$, and $V_2$ in the same biliaison class as $V_1$, where the
map ${\mathcal E}''(-b_2) \rightarrow {\mathcal N}$ extends to ${\mathcal
E}(-b_2)$.  We may assume $b_2 \gg b_1$ and then change the map from ${\mathcal
L}(-b_1) \rightarrow {\mathcal  N}$ again so that $b_1 = b_2 = b$.

Consider the  diagram
\setcounter{equation}{0}
\begin{equation}
\begin{array}{ccccccccc}
& &0 & &0 \\
& &\downarrow & &\downarrow \\
& &{\mathcal E}'' &= &{\mathcal E}'' \\
& &\downarrow & &\downarrow \\
0 &\rightarrow &{\mathcal E}' &\rightarrow &{\mathcal E} &\rightarrow &{\mathcal
F} &\rightarrow &0 \\
& &\downarrow & &\downarrow & &\| \\
0 &\rightarrow &{\mathcal L} &\rightarrow &{\mathcal R} &\rightarrow &{\mathcal
F} &\rightarrow &0 \\
& &\downarrow & &\downarrow \\
& &0 & &0.
\end{array}
\end{equation}

Consider also the diagram
\begin{equation}
\begin{array}{ccccc}
& &0 & &0 \\
& &\downarrow & &\downarrow \\
& &{\mathcal E}_0 \oplus {\mathcal E}''(-b) &= &{\mathcal E}_0 \oplus {\mathcal
E}''(-b) \\
& &\downarrow & &\downarrow \\
0 &\rightarrow &{\mathcal E}_0 \oplus {\mathcal E}'(-b) &\rightarrow &{\mathcal
N} \\
& &\downarrow & &\downarrow \\
& &{\mathcal L}(-b) &\rightarrow &{\mathcal G} \\
& &\downarrow & &\downarrow \\
& &0 & &0
\end{array}
\end{equation}
using the map ${\mathcal E}'(-b)$ to ${\mathcal N}$ obtained by restricting  from
${\mathcal E}(-b)$ to ${\mathcal N}$.  We obtain a map $s_0: {\mathcal L}(-b)
\rightarrow {\mathcal G}$ on the quotients with no special property.  (Note rank
${\mathcal G} = 2$.)

Taking ${\mathcal H}om(\cdot,{\mathcal N}(b))$ of the first diagram, and
composing with the map ${\mathcal N} \rightarrow {\mathcal G}$ from the second
diagram, we get
\[
\begin{array}{ccc}
{\mathcal H}om({\mathcal R},{\mathcal N}(b)) &\rightarrow &{\mathcal
H}om({\mathcal E},{\mathcal N}(b)) \\
\downarrow & &\downarrow \\
{\mathcal H}om({\mathcal L},{\mathcal N}(b)) &\rightarrow &{\mathcal
H}om({\mathcal E}',{\mathcal N}(b)) \\
\downarrow \\
{\mathcal H}om({\mathcal L},{\mathcal G}(b)).
\end{array}
\]
We may assume that $b$ is so large that ${\mathcal H}om({\mathcal R},{\mathcal
N}(b))$ is generated by global sections.

Let $W \subseteq H^0({\mathcal H}om({\mathcal L},{\mathcal G}(b))$ be the image
of $H^0({\mathcal H}om({\mathcal R},{\mathcal N}(b))$ above.  As before, since
the sheaves ${\mathcal R}$ and ${\mathcal L}$ are locally free in codimension
$1$, the cokernels of these maps of sheaves are ${\mathcal E}xt$ sheaves with
support in codimension $\ge 2$, so $W$ will satisfy the conditions of $(3.2)$ and
generate a subsheaf of rank $2$.  Therefore a general element $s \in W$ will give
a map ${\mathcal L}(-b) \rightarrow {\mathcal G}$ whose quotient is torsion-free
and locally free in codimension $1$.  However, instead of using this map, we add
it to the existing map $s_0$ of ${\mathcal L}(-b)$ to ${\mathcal G}$, and the
proof of $(3.2)$ shows that for general $s \in W$, the sum $s_0 + s$ will give a
good cokernel.  Since by construction $s$ comes from a map of ${\mathcal R}$ to
${\mathcal N}(b)$, we get maps of ${\mathcal E}(-b)$ to ${\mathcal N}$ and
${\mathcal E}'(-b)$ to ${\mathcal  N}$, which we add to the existing maps.  Thus
we get a diagram like $(2)$ above, but with new maps, where the cokernel of the
bottom two rows is
${\mathcal I}_{V'}(a')$, with
$V'$ in the same biliaison class as $V_2$, for some $a' \in {\mathbb Z}$, and
where the map ${\mathcal E}'(-b)$ to
${\mathcal N}$ extends to ${\mathcal E}(-b)$, as required.

\bigskip
\noindent
{\bf Lemma 4.6.} {\em Again with $X$ normal {\em ACM}, suppose given $V$ of
codimension $2$ and exact sequences}
\[
0 \rightarrow {\mathcal E}_0 \oplus {\mathcal E}' \oplus {\mathcal L} \rightarrow
{\mathcal N} \rightarrow {\mathcal I}_V(a) \rightarrow 0
\]
\[
0 \rightarrow {\mathcal E}' \rightarrow {\mathcal E} \rightarrow {\mathcal L}
\rightarrow 0
\]
{\em with ${\mathcal E}_0,{\mathcal N}$ coherent, ${\mathcal E}',{\mathcal E}$
layered {\em ACM} sheaves and ${\mathcal L}$ a {\em rank} $1$ {\em ACM} sheaf. 
Then for any $b \gg 0$ there is an exact sequence}
\[
0 \rightarrow {\mathcal E}_0 \oplus {\mathcal E}(-b) \rightarrow {\mathcal N}
\rightarrow {\mathcal I}_{V'}(a') \rightarrow 0
\]
{\em with $V'$ in the same biliaison class as $V$.}

\bigskip
\noindent
{\em Proof.}  By $(4.5)$ we can find
\[
0 \rightarrow {\mathcal E}_0 \oplus {\mathcal E}'(-b_1) \oplus {\mathcal L}
\rightarrow {\mathcal N} \rightarrow {\mathcal I}_{V_1}(a_1) \rightarrow 0
\]
where the map ${\mathcal E}'(-b_1)$ to ${\mathcal N}$ extends to a map of
${\mathcal E}(-b_1)$ to ${\mathcal N}$.  Next, by the method of $(3.4)$ we
replace ${\mathcal L}$ by ${\mathcal L}(-b_2)$, and we may assume $b_1 = b_2 = b$
by taking both sufficiently large.  Now the same proof as for the latter part of
$(4.5)$ shows how to get the desired sequence.

\bigskip
\noindent
{\em Proof of $(4.3)$.}  Let $V$ be a codimension $2$ ACM subscheme of $X$.  Take
a set of generators for $H_*^0({\mathcal I}_V)$ and thus obtain an ${\mathcal
E}$-type resolution
\[
0 \rightarrow {\mathcal E} \rightarrow {\mathcal L} \rightarrow {\mathcal I}_V
\rightarrow 0
\]
with ${\mathcal L}$ dissoci\'e and $H_*^1({\mathcal E}) = 0$.  Since $V$ is ACM,
it follows that ${\mathcal E}$ is an ACM sheaf on $X$.

By $(4.4)$, using condition (B) and hence (B${}^*$), we can find a resolution
\[
0 \rightarrow {\mathcal E} \rightarrow {\mathcal G}_1 \rightarrow {\mathcal G}_2
\rightarrow 0
\]
with ${\mathcal G}_1,{\mathcal G}_2$ layered ACM sheaves.  Taking ${\mathcal H}$
to be the fibered sum of ${\mathcal L}$ and ${\mathcal G}_1$ over ${\mathcal E}$, 
we obtain sequences 
\[
0 \rightarrow {\mathcal G}_1 \rightarrow {\mathcal H} \rightarrow {\mathcal I}_V
\rightarrow 0
\]
and
\[
0 \rightarrow {\mathcal L} \rightarrow {\mathcal H} \rightarrow {\mathcal G}_2
\rightarrow 0.
\]
Thus ${\mathcal H}$ is layered, and we have a resolution of ${\mathcal I}_V$ with
${\mathcal G}_1,{\mathcal H}$ both  layered ACM sheaves.  (In fact, since
${\mathcal L}$ is dissoci\'e and ${\mathcal G}_2$ is ACM, this last sequence
splits, and ${\mathcal H} \cong {\mathcal L} \oplus {\mathcal G}_2$, but we do
not need this.)

Let $V'$ be another codimension $2$ ACM subscheme of $X$, and let
\[
0 \rightarrow {\mathcal G}'_1 \rightarrow {\mathcal H}' \rightarrow {\mathcal
I}_{V'} \rightarrow 0
\]
be a similar sequence for $V'$.  Now adding ${\mathcal H}'$ to the sequence for
$V$ and ${\mathcal H}$ to the sequence for $V'$ we obtain resolutions
\[
0 \rightarrow {\mathcal G}_1 \oplus {\mathcal H}' \rightarrow {\mathcal H} \oplus
{\mathcal H}' \rightarrow {\mathcal I}_V \rightarrow 0
\]
\[
0 \rightarrow {\mathcal G}'_1 \oplus {\mathcal H} \rightarrow {\mathcal H} \oplus
{\mathcal H}' \rightarrow {\mathcal I}_{V'} \rightarrow 0
\]
by layered ACM sheaves, with the same sheaf in the middle of each.  Under
hypothesis a) we add $({\mathcal H} \oplus {\mathcal H}')^{\vee}$ to both
sequences, or under hypothesis b) we add a suitable sum $\oplus {\mathcal L}_i$
of rank $1$ ACM sheaves to each, and then we obtain sequences (changing notation)
\[
0 \rightarrow {\mathcal E} \rightarrow {\mathcal N} \rightarrow {\mathcal I}_V
\rightarrow 0
\]
\[
0 \rightarrow {\mathcal E}' \rightarrow {\mathcal N} \rightarrow {\mathcal
I}_{V'} \rightarrow 0
\]
where ${\mathcal E},{\mathcal E}'$ and ${\mathcal N}$ are all three orientable
layered ACM sheaves on $X$.

Next, we invoke condition (C) to tell us that ${\mathcal E}$ and ${\mathcal E}'$
are both equal to $r \cdot {\mathcal O}$ in the group $G'$, and hence equal to
each other in $G'$.  What does this equality mean?  It means that up to twist, we
can transform ${\mathcal E}$ to ${\mathcal E}'$ by a finite number of operations
of adding and subtracting expressions $({\mathcal F} - {\mathcal F}' - {\mathcal
F}'')$ whenever $0 \rightarrow {\mathcal F}' \rightarrow {\mathcal F} \rightarrow
{\mathcal F}'' \rightarrow 0$ is an exact sequence in the category ${\mathcal M}$
of layered ACM sheaves on $X$.  In fact, it is enough to consider such
expressions in which ${\mathcal F}''$ has rank $1$, because any exact sequence $0
\rightarrow {\mathcal F}' \rightarrow {\mathcal F} \rightarrow {\mathcal F}''
\rightarrow 0$ of layered ACM sheaves can be resolved into a finite number of
such sequences whose last term has rank $1$, by using the filtration on
${\mathcal F}''$.  By moving negative terms to the other side of the equation, we
find a layered ACM sheaf ${\mathcal G}$ such that one can transform ${\mathcal E}
\oplus {\mathcal G}$ into ${\mathcal E}' \oplus {\mathcal G}$ by a finite number
of substitutions of ${\mathcal F}$ for ${\mathcal F}' \oplus {\mathcal F}''$ or
vice versa, when $0 \rightarrow {\mathcal F}' \rightarrow {\mathcal F}
\rightarrow {\mathcal F}'' \rightarrow 0$ is an exact sequence as above.  All
this takes place up to twists.

Now $(3.4)$ allows us to replace a layered ACM sheaf by the direct sum of its
factors (sufficiently twisted), and $(4.6)$ allows us to replace a direct sum
${\mathcal E}' \oplus {\mathcal L}$ by an extension ${\mathcal E}$ of ${\mathcal
L}$ by ${\mathcal E}'$, sufficiently twisted.  By applying these lemmas to the
resolution
\[
0 \rightarrow {\mathcal E} \oplus {\mathcal G} \rightarrow {\mathcal
N} \oplus {\mathcal G} 
\rightarrow {\mathcal I}_V
\rightarrow 0
\]
of $V$ we thus obtain a resolution
\[
0 \rightarrow {\mathcal E}'(-b) \oplus {\mathcal G}_1 \rightarrow {\mathcal N}
\oplus {\mathcal G} \rightarrow {\mathcal I}_{V_1}(a_1) \rightarrow 0
\]
for some $V_1$ in the biliaison class of $V$, and where ${\mathcal G}$ is a
layered ACM sheaf, and ${\mathcal G}_1$ differs from ${\mathcal G}$ only by
twists of its direct summands.

We compare this to the resolution
\[
0 \rightarrow {\mathcal E}' \oplus {\mathcal G} \rightarrow {\mathcal N} \oplus
{\mathcal G} \rightarrow {\mathcal I}_{V'} \rightarrow 0
\]
of $V'$.  Note that the sheaves ${\mathcal E}'(-b) \oplus {\mathcal G}_1$ and
${\mathcal E}' \oplus {\mathcal G}$ have all the same rank $1$ factors up to
twist.  So from $(3.1)$ we
conclude that $V$, $V_1$ and $V'$ are all in the same biliaison class.

\bigskip
\noindent
{\bf Theorem 4.7.}  {\em Let $X$ be a three-dimensional normal arithmetically
Gorenstein projective scheme.  Assume condition {\em (C)} and the stronger
condition {\em (B$''$)}:  every orientable {\em ACM} sheaf on $X$ is layered. 
Then two locally Cohen--Macaulay curves $C_1$, $C_2$ in $X$ are in the same
Gorenstein biliaison class if and only if their Rao modules $H_*^1({\mathcal
I}_{C_1})$ and $H_*^1({\mathcal I}_{C_2})$ are isomorphic, up to twist.}

\bigskip
\noindent
{\em Proof.}  One direction is elementary, since we do biliaison on ACM surfaces
$(1.2)$.

For the other direction, let $C$ be any locally Cohen--Macaulay curve on $X$,
with Rao module $M = H_*^1({\mathcal I}_C)$.  Take an ${\mathcal N}$-type 
resolution of $C$, that is, an exact sequence
\[
0 \rightarrow {\mathcal L} \rightarrow {\mathcal N} \rightarrow {\mathcal I}_C
\rightarrow 0
\]
with ${\mathcal L}$ dissoci\'e, and ${\mathcal N}$ a locally Cohen--Macaulay
sheaf with $H_*^1({\mathcal N}) = M$ and $H_*^2({\mathcal N}) = 0$
\cite[1.12,1.13]{RTLRP}.  Let
$N = H_*^0({\mathcal N})$.  As in \cite[3.2]{RTLRP} a two-step resolution of $N$
\[
0 \rightarrow P \rightarrow L_1 \rightarrow L_0 \rightarrow N \rightarrow 0
\]
over the homogeneous coordinate ring $S$ of $X$ gives a maximal Cohen--Macaulay
module $P$, and dualizing this sequence gives a sequence
\[
0 \rightarrow N^{\vee} \rightarrow L_0^{\vee} \rightarrow L_1^{\vee} \rightarrow
P^{\vee}
\rightarrow M^* \rightarrow 0,
\]
where $\cdot^{\vee} = \mbox{Hom}_S(\cdot,S)$ and $M^* = \mbox{Hom}_k(M,k)$.

Now take a minimal free resolution of $M^*$,
\[
0 \rightarrow R \rightarrow F_2 \rightarrow F_1 \rightarrow F_0 \rightarrow M^*
\rightarrow 0
\]
defining a new $S$-module $R$.  Let ${\mathcal N}'$ be the sheaf associated to
$R^{\vee}$.  There is a natural map of the free resolution of $M^*$ into the
earlier resolution, and so we get a map $R \rightarrow N^{\vee}$ giving rise to a
map of sheaves ${\mathcal N} \rightarrow {\mathcal N}'$, which induces an
isomorphism on the $H_*^1$-modules, both isomorphic to $M$.  By adding suitable
dissoci\'e sheaves to the original ${\mathcal L}$ and ${\mathcal N}$, we may
assume that
${\mathcal N} \rightarrow {\mathcal N}'$ is surjective also on $H_*^0$, so we get
\[
0 \rightarrow {\mathcal E} \rightarrow {\mathcal N} \rightarrow {\mathcal N}'
\rightarrow 0
\]
with ${\mathcal E}$ an ACM sheaf on $X$.

Let ${\mathcal L}'$ be a dissoci\'e sheaf mapping to ${\mathcal N}'$
\[
0 \rightarrow {\mathcal L}' \rightarrow {\mathcal N}' \rightarrow {\mathcal
I}_{C'}(a') \rightarrow 0
\]
with cokernel ${\mathcal I}_{C'}(a')$ for some curve $C'$.  Composing with the
map ${\mathcal N} \rightarrow {\mathcal N}'$, we get
\[
0 \rightarrow {\mathcal E}' \rightarrow {\mathcal N} \rightarrow {\mathcal
I}_{C'}(a') \rightarrow 0
\]
where ${\mathcal E}'$ is an extension of ${\mathcal L}'$ by ${\mathcal E}$
(actually a direct sum), hence another ACM sheaf.

Now by hypothesis (B$''$), ${\mathcal E}'$ is a layered ACM sheaf, and it is also
orientable, since ${\mathcal N}$ is.

At this point we can repeat the last part of the proof of $(4.3)$, using
condition (C), and noting that for that part of the proof, it was not necessary
for ${\mathcal N}$ to be ACM.  Thus $C$ and $C'$ are in the same Gorenstein
biliaison class.  Since the biliaison class of $C'$ depends only on the Rao
module, and not on $C$, this shows that any two curves with the same Rao module
are in the same biliaison class.

\section{Applications I: Schemes of finite representation type}
\label{sec5}

We say a projective ACM scheme $X$ has {\em finite representation type} if there
are only finitely many isomorphism classes of indecomposable ACM sheaves on $X$
up to twist.  The projective schemes of finite representation type have been
classified by Eisenbud and Herzog, see \cite[17.10]{Y}.  Those of dimension $\ge
2$ are (i) ${\mathbb P}^n$ for $n \ge 2$, (ii) a non-singular quadric hypersurface
$Q^n$ in ${\mathbb P}^{n+1}$ for $n \ge 2$, (iii) the rational cubic scroll in
${\mathbb P}^4$, and (iv) the Veronese surface in ${\mathbb P}^5$.

If $X$ is ${\mathbb P}^n$, our conditions (B) and (C) are trivial, because the
only ACM sheaves on $X$ are dissoci\'e.  Theorem $(4.3)$ tells us that all
codimension $2$ ACM schemes are in the same Gorenstein biliaison class, and
$(4.7)$ tells us in ${\mathbb P}^3$ that the Gorenstein biliaison class of any
curve is determined by its Rao module.  Since Gorenstein biliaison in codimension
$2$ in ${\mathbb P}^n$ is just complete intersection biliaison, we recover the
well-known results of Ap\'ery--Gaeta--Peskine--Szpiro \cite{PS} and of Rao
\cite{R}.  Our proof cannot be considered a new proof, because the methods we use
are generalizations of the earlier proofs.

If $X$ is a non-singular quadric hypersurface $Q^n$ in ${\mathbb P}^{n+1}$, over
an algebraically closed field $k$ of characteristic $\ne 2$, the non-free
indecomposable ACM sheaves have been classified by the work of Buchweitz,
Eisenbud, and Herzog \cite{BEH} and Kn\"orrer \cite{K}, see \cite[14.10]{Y}.  Up
to isomorphism and twist there is just one non-free ACM sheaf if $n$ is odd, or
two if $n$ is even.  The rank of these sheaves is $2^m$ where $m = \left[ \frac
{n-1}{2} \right]$.

\bigskip
\noindent
{\bf Proposition 5.1.}  {\em If $X$ is a non-singular quadric surface in
${\mathbb P}^3$, all zero-dimensional subschemes of $X$ are in the same
Gorenstein biliaison class.}

\bigskip
\noindent
{\em Proof.}  In this case there are two non-trivial 
indecomposable ACM sheaves of rank $1$, namely  ${\mathcal I}_L$ and ${\mathcal
I}_M$ where $L,M$ represent the two classes of lines on $X$.  Thus any ACM sheaf
is a direct sum of rank $1$ ACM sheaves, and condition (B) is trivial.

For condition (C), note that $L+M=H$, the hyperplane class, and since two lines
in the same family do not meet, there is an exact sequence
\[
0 \rightarrow {\mathcal I}_L \rightarrow {\mathcal O}^2 \rightarrow {\mathcal
I}_M(1) \rightarrow 0.
\]
If ${\mathcal E}$ is an orientable ACM sheaf, it must have the same number of
$L$'s and $M$'s in its direct sum decomposition, so the existence of this
sequence proves (C).

We conclude the result from $(4.3)$.

\bigskip
\noindent
{\bf Remark 5.2.}  This generalizes an earlier result \cite[2.3]{SEG} which
showed that any set of $n$ general points on the quadric surface was in the
biliaison class of a point.

\bigskip
\noindent
{\bf Example 5.3.}  If $X$ is a non-singular quadric hypersurface of dimension $n
\ge 3$ in ${\mathbb P}^{n+1}$, then $\mbox{Pic } X = {\mathbb Z}$, generated by
${\mathcal O}(1)$, and there are no non-trivial extensions, so condition (C)
holds trivially.  On the other hand, there are indecomposable ACM sheaves of rank
$\ge 2$.  Any sequence as in condition (B) would have to split, so condition (B)
fails.  We conclude there are many inequivalent Gorenstein biliaison classes of
codimension $2$ subschemes of $X$.  In fact, since $\mbox{Pic } X = {\mathbb Z}$,
Gorenstein biliaison is the same as CI-biliaison, and the biliaison classes are
determined by stable equivalence classes (up to twist) of ACM sheaves on $X$
\cite[2.4]{RTLRP}.

\bigskip
\noindent
{\bf Proposition 5.4.}  {\em All zero-dimensional subschemes of the rational
cubic scroll $X$ in ${\mathbb P}^4$ are in the same Gorenstein biliaison class.}

\bigskip
\noindent
{\em Proof.}  First recall that $\mbox{Pic } X \cong {\mathbb Z}^2$ and is
generated by the hyperplane class $H$ and the class of a fiber $F$ of the ruled
surface.  It is well known \cite[5.10]{N} that up to twist the only rank $1$ ACM
sheaves are ${\mathcal O}_X$, ${\mathcal I}_F$, ${\mathcal I}_{H-F}$, ${\mathcal
I}_{H-2F}$.  Here $H-F$ is a conic, and $H-2F$ is the exceptional line $E$ with
$E^2 = -1$.

Looking at the proof \cite[16.12]{Y} that $X$ is of finite representation type,
we find that there is just one (up to twist) indecomposable ACM sheaf of rank $>
1$, it has rank $2$, and is obtained as the first syzygy of ${\mathcal I}_F$:
\[
0 \rightarrow {\mathcal E}_0 \rightarrow {\mathcal O}(-1)^3 \rightarrow {\mathcal
I}_F \rightarrow 0.
\]
Note that $c_1({\mathcal E}_0) = F - 3H$.

First we will show that ${\mathcal E}_0$ is isomorphic to an extension of rank
$1$ ACM sheaves
\setcounter{equation}{0}
\begin{equation}
0 \rightarrow {\mathcal I}_F(-1) \rightarrow {\mathcal F} \rightarrow {\mathcal
I}_{H-2F}(-1) \rightarrow 0.
\end{equation}
We compute $\mbox{Ext}^2({\mathcal I}_{H-2F}(-1),{\mathcal I}_F(-1)) =
H^1({\mathcal I}_{3F}(1))$ is of dimension $1$, so there exists a non-trivial
extension $(1)$ and it is an ACM sheaf of rank $2$.  Next, note that
$h^0({\mathcal F}(1)) = 0$, $h^0({\mathcal F}(2)) = 6$.  So if ${\mathcal F}$
were decomposable, then ${\mathcal F} = {\mathcal L}_1 \oplus {\mathcal L}_2$,
where ${\mathcal L}_1,{\mathcal L}_2$ are rank $1$ ACM sheaves with
$h^0({\mathcal L}_i(1)) = 0$, and at least one of $h^0({\mathcal L}_i(2)) \ne
0$.  The only rank $1$ ACM sheaves ${\mathcal L}$ with $h^0({\mathcal L}(1)) =
0$, $h^0({\mathcal L}(2)) \ne 0$ are ${\mathcal O}(-2)$, ${\mathcal I}_F(-1)$,
${\mathcal I}_{H-F}(-1)$, ${\mathcal I}_{H-2F}(-1)$.  For these sheaves
${\mathcal L}$ we have $h^0({\mathcal L}(2)) = 1,3,2,3$, respectively.  Taking
into account $c_1({\mathcal F}) = F-3H$, the only possibility would be ${\mathcal
F} \cong {\mathcal I}_F(-1) \oplus {\mathcal I}_{H-2F}(-1)$, which is impossible
since the sequence $(1)$ is non-split.

Therefore ${\mathcal F}$ is an indecomposable rank $2$ ACM sheaf on $X$.  Since
there is only one such up to twist \cite[16.12]{Y}, we find checking Chern
classes that ${\mathcal F} \cong {\mathcal E}_0$.  Thus ${\mathcal E}_0$ is
already layered, and condition (B) of $(4.3)$ is satisfied.

Note that $X$ is not arithmetically Gorenstein, but that the rank $1$ ACM
divisors $H$, $F$, $H-F$ and their twists generate $\mbox{Pic } X$ as a monoid,
so condition b) of $(4.3)$ is satisfied.

For condition (C), note that since $F \cdot F = 0$, there is an exact sequence
\begin{equation}
0 \rightarrow {\mathcal I}_F \rightarrow {\mathcal O}^2 \rightarrow {\mathcal
I}_{H-F}(1) \rightarrow 0,
\end{equation}
and since $F \cdot 2F = 0$ another sequence
\begin{equation}
0 \rightarrow {\mathcal I}_F \rightarrow {\mathcal O}_X \oplus {\mathcal
I}_{H-F}(1) \rightarrow {\mathcal I}_{H-2F}(1) \rightarrow 0.
\end{equation}

Now suppose given an orientable ACM sheaf on $X$.  Up to twists we can write it
as a direct sum
\[
a \cdot {\mathcal I}_F + b \cdot {\mathcal I}_{H-F} + c{\mathcal I}_{H-2F} + d
\cdot {\mathcal E}_0 + e \cdot {\mathcal O}
\]
with $a,b,c,d,e \ge 0$.

In the quotient Grothendieck group $G'$, we use sequence $(1)$ to replace
${\mathcal E}_0$ by ${\mathcal I}_F + {\mathcal I}_{H-2F}$ and then we may assume
$d = 0$.  Since the sheaf is orientable, if $c \ne 0$, then also $a \ne 0$.  So
we can use sequence $(3)$ to replace occurrences of ${\mathcal I}_{H-2F} +
{\mathcal I}_F$ by ${\mathcal I}_{H-F} + {\mathcal O}$ and then assume
$c=0$.  Now if $b$ is non-zero, it must be equal to $a$, and we use sequence
$(2)$ to reduce to a multiple of ${\mathcal O}$.  Thus condition (C) is satisfied.

Thus by $(4.3)$ every zero-scheme on $X$ is in the biliaison class of a point.

\bigskip
\noindent
{\bf Remark 5.5.}  This generalizes an earlier result \cite[3.4]{C} that treated
only sets of points in general position.

\bigskip
\noindent
{\bf Proposition 5.6.} {\em On the Veronese surface $X$ in ${\mathbb P}^5$ there
are infinitely many biliaison classes of zero-schemes, indexed by the even
integers.}

\bigskip
\noindent
{\em Proof.}  The Veronese surface is the $2$-tuple embedding of ${\mathbb P}^2$
in ${\mathbb P}^5$.  Thus $\mbox{Pic } X = {\mathbb Z}$, generated by the image
of a line of ${\mathbb P}^2$, which becomes a conic $C \subseteq X$.  The
hyperplane class is $H = 2C$.  Every curve on $X$ is an ACM curve.  $X$ is not
arithmetically Gorenstein, but its Picard group is generated as a monoid by $C$
and $H-C$, which are ACM curves.  There are no non-trivial extensions of rank $1$
ACM sheaves, so property (C) fails, since ${\mathcal I}_C \oplus {\mathcal I}_C$
is an orientable ACM sheaf not equivalent in the group $G'$ to $2 \cdot {\mathcal
O}$.

Following the proof of \cite[16.10]{Y} we see that there is just one (up to
twist) indecomposable ACM sheaf of rank $> 1$; it has rank $2$ and is obtained as
the first syzygy of ${\mathcal I}_C$:
\[
0 \rightarrow {\mathcal E}_0 \rightarrow {\mathcal O}(-1)^3 \rightarrow {\mathcal
I}_C \rightarrow 0.
\]
Pulling this sequence back to ${\mathbb P}^2$ one recognizes that ${\mathcal E}_0
\cong \Omega_{{\mathbb P}^2}^1(-1)$.  Note that $H^1({\mathbb
P}^2,\Omega_{{\mathbb P}^2}^1) \ne 0$, but still ${\mathcal E}_0$ is an ACM sheaf
on $X$ because twists on $X$ correspond to even twists on ${\mathbb P}^2$, and no
other twist of $\Omega_{{\mathbb P}^2}^1$ has a non-zero $H^1$.  Now ${\mathcal
E}_0(2)$ is the tangent bundle on $X$, and transporting a well-known sequence
with the tangent bundle from ${\mathbb P}^2$ we find
\[
0 \rightarrow {\mathcal O}(-2) \rightarrow {\mathcal I}_C(-1)^3 \rightarrow
{\mathcal E}_0 \rightarrow 0,
\]
so condition (B) is satisfied.

Since condition (C) fails, $(4.2)$ tells us that there is more than one biliaison
class of zero-schemes on $X$.  This is clear anyway, because $H.Y$ is even
for any curve $Y$ on $X$, so biliaison preserves the parity of the length of a
zero-scheme $Z$.

To investigate more exactly what the biliaison classes are, let $Z,Z'$ be any two
zero-schemes in $X$.  Since condition (B) holds, we can use $(4.4)$ and the
beginning of the proof of $(4.3)$ to find resolutions
\[
0 \rightarrow {\mathcal E} \rightarrow {\mathcal N} \rightarrow {\mathcal I}_Z
\rightarrow 0
\]
\[
0 \rightarrow {\mathcal E}' \rightarrow {\mathcal N}' \rightarrow {\mathcal
I}_{Z'} \rightarrow 0
\]
where ${\mathcal E},{\mathcal E}'$ ${\mathcal N},{\mathcal N}'$ are orientable
layered ACM sheaves, and we may even assume ${\mathcal N} = {\mathcal N}'$.  In
our case, this means ${\mathcal E},{\mathcal E}',{\mathcal N},{\mathcal N}'$ are
each direct sums of copies of ${\mathcal O}(a_i)$ and ${\mathcal I}_C(b_i)$ for
various $a_i,b_i$.

Without assuming ${\mathcal N} = {\mathcal N}'$, if $Z$ and $Z'$ are in the same
biliaison class, then the exact sequences of $(4.1)$ show that there are layered
ACM sheaves ${\mathcal F},{\mathcal F}'$ with the same rank $1$ factors, up to
twist, and an equivalence ${\mathcal N} + {\mathcal E}' + {\mathcal F} =
{\mathcal N}' + {\mathcal E} + {\mathcal F}'$ in the Grothendieck group $G'$. 
Since ${\mathcal F}$ and ${\mathcal F}'$ have the same factors, up to twist, we
find ${\mathcal N} + {\mathcal E}' = {\mathcal N}' + {\mathcal E}$ in $G'$.

Let $m(Z) =$ the number of copies of ${\mathcal I}_C$ in the direct sum
decomposition of ${\mathcal N}$, minus the number in ${\mathcal E}$.  Since
${\mathcal E},{\mathcal N}$ are orientable, $m(Z)$ is an even integer.  What
$(4.1)$ tells us in this case is that $m(Z)$ is an invariant of the biliaison
class.

Conversely, if $Z,Z'$ are two schemes with $m(Z) = m(Z')$, then as in the proof
of $(4.3)$, we find the sequences above with ${\mathcal  N} = {\mathcal N}'$. 
Then ${\mathcal E}$ and ${\mathcal E}'$ have the same number of direct summands
${\mathcal I}_C(a_i)$ for various $a_i$, and it follows from $(3.1)$ that they
are in the same biliaison class.

Given $m > 0$, $m$ even, take ${\mathcal N} = \oplus_{i=1}^m {\mathcal I}_C$,
choose a dissoci\'e sheaf ${\mathcal L}$ such that
\[
0 \rightarrow {\mathcal L} \rightarrow {\mathcal N} \rightarrow {\mathcal I}_Z(a)
\rightarrow 0
\]
for some zero-scheme $Z$, and then $Z$ will have invariant $m$.

If $m < 0$, $m$ even, let $n = -m$.  Take a general $n \times (n+1)$ matrix of
linear forms on ${\mathbb P}^2$, and let $Z$ be the associated determinantal
scheme, so that we have a resolution
\[
0 \rightarrow {\mathcal O}_{{\mathbb P}^2}(-n-1)^n \rightarrow {\mathcal
O}_{{\mathbb P}^2}(-n)^{n+1} \rightarrow {\mathcal I}_{Z,{\mathbb P}^2}
\rightarrow 0.
\]
Transporting this to the Veronese surface gives $Z$ with resolution
\[
0 \rightarrow {\mathcal I}_C\left( -\frac {n}{2}\right)^n \rightarrow {\mathcal
O}\left( -\frac {n}{2}\right)^{n+1} \rightarrow {\mathcal I}_Z \rightarrow 0
\]
on the Veronese surface, showing that $Z$ has invariant $-n = m$.

Thus the biliaison classes of zero-schemes on $X$ are in one-to-one
correspondence with the even integers via the invariant $m(Z)$.

\bigskip
\noindent
{\bf Example 5.7.}  Let $X$ be the Veronese surface in ${\mathbb P}^5$.  If $Z =
P$ is a single point, there is a resolution
\[
0 \rightarrow {\mathcal O}(-1) \rightarrow {\mathcal I}_C \oplus {\mathcal I}_C
\rightarrow {\mathcal I}_P \rightarrow 0,
\]
since two lines in ${\mathbb P}^2$ meet in a point.  Thus $m(P) = 2$.  It is easy
to see that
$Z =$ three points in general position has $m(Z) = -2$.  On the other hand, the
image in $X$ of $3$ points on a line in ${\mathbb P}^2$ is in the biliaison class
of a point, so $m = 2$.

If $Z$ is two points, there is a resolution
\[
0 \rightarrow {\mathcal I}_C(-1) \rightarrow {\mathcal O}(-1) \oplus {\mathcal
I}_C \rightarrow {\mathcal I}_Z \rightarrow 0
\]
since a line and a conic in ${\mathbb P}^2$ meet in $2$ points.  Hence $m(Z) =
0$, and we conclude that $2$ points are in the biliaison class of a complete
intersection.  This is not obvious!  We leave as an amusing exercise for the
reader to find explicit biliaisons (among non-empty zero schemes) that relate $2$
points to a complete intersection in $X$.

To find a scheme with $m(Z) = 4$, we use a sequence
\[
0 \rightarrow {\mathcal O}(-1)^3 \rightarrow {\mathcal I}_C^4 \rightarrow
{\mathcal I}_Z(1) \rightarrow 0.
\]
Then $Z$ is a set of $6$ points, which may be taken in general position.

\section{Applications, II: Quadric cones}
\label{sec6}

In this section we consider some schemes that are not of finite representation
type, but on which we can still determine all possible ACM sheaves.  We use the
matrix factorizations of Eisenbud \cite{E}, see \cite[Ch.~7]{Y} and the
periodicity theorems of Kn\"orrer \cite{K}, see \cite[Ch.~12]{Y}.

\bigskip
\noindent
{\bf Proposition 6.1.}  {\em Let $X$ be a quadric cone in ${\mathbb P}^3$.  Then
all zero-schemes on $X$ are in the same Gorenstein biliaison class.}

\bigskip
\noindent
{\em Proof.}  Recall that $\mbox{APic } X = {\mathbb Z}$, generated by a line
$L$, and $2L = H$ is the hyperplane class.  The surface is arithmetically
Gorenstein, and there is an exact sequence
\[
0 \rightarrow {\mathcal I}_L \rightarrow {\mathcal O}^2 \rightarrow {\mathcal
I}_L(1) \rightarrow 0,
\]
so condition (C) is satisfied.

To determine ACM sheaves of higher rank, we use a result of Buchweitz, Greuel,
and Schreyer \cite[4.1]{BGS}.  If $S = k[x,t]$ and $R = S/(x^2)$, they found that
there are three types of indecomposable maximal Cohen--Macaulay modules on $R$,
given by matrix factorizations $(\varphi,\psi)$

\begin{itemize}
\item[(i)] $((x^2),(1))$
\item[(ii)] $((x),(x))$ 
\item[(iii)] $(\varphi_{\ell},\psi_{\ell})$ where $\varphi_{\ell} = \psi_{\ell} =
\begin{pmatrix}
x & t^{\ell} \\
0 & -x
\end{pmatrix}$ for $\ell = 1,2,3,\dots$.
\end{itemize}
The corresponding $R$-module will be the cokernel of a map of free $S$-modules
determined by the first matrix in each pair.  Thus (i) gives $S/(x^2) = R$, (ii)
gives $S/(x) = R/(x)$, and (iii) gives a module isomorphic to the ideal
$(x,t^{\ell})$ in $R$.

By Kn\"orrer's periodicity theorem \cite[3.1]{K}, if we let $S_2 = k[x,u,v,t]$
and $R_2 = S_2/(x^2+uv)$, then the MCM modules on $R'$ are given by matrix
factorizations
\[
\left( \begin{pmatrix}
u & \psi\\
\varphi & -v
\end{pmatrix}, \begin{pmatrix}
v & \psi\\
\varphi & -u
\end{pmatrix} \right)
\]
for each matrix factorization $(\varphi,\psi)$ of an MCM module on $R$.

Under this correspondence, type (i) gives us $R_2$; type (ii) gives us the ideal
$I_L = (x,u)$ of a line in $X$, and type (iii) gives a module $M_{\ell}$, for
$\ell = 1,2,3,\dots$, given by a matrix factorization
\[
\left( \begin{pmatrix}
u & 0 & x & t^{\ell} \\
0 & u & 0 & -x \\
x & t^{\ell} & -v & 0 \\
0 & -x & 0 & -v
\end{pmatrix}, \begin{pmatrix}
v & 0 & x & t^{\ell} \\
0 & v & 0 & -x \\
x & t^{\ell} & -u & 0 \\
0 & -x & 0 & -u
\end{pmatrix} \right)\ .
\]
This is equivalent, permuting columns and rows to
\[
\left( \begin{pmatrix}
u & x & 0 & t^{\ell} \\
x & -v & t^{\ell} & 0 \\
0 & 0 & u & -x \\
0 & 0 & -x & -v
\end{pmatrix}, \begin{pmatrix}
v & x & 0 & t^{\ell} \\
x & -u & t^{\ell} & 0 \\
0 & 0 & v & -x \\
0 & 0 & -x & u
\end{pmatrix}\right)\ .
\]
Thus we see that $M_{\ell}$ is an extension
\[
0 \rightarrow I_L \rightarrow M_{\ell} \rightarrow I_L(-\ell+1) \rightarrow 0
\]
of rank $1$ $R_2$-modules.  The corresponding ACM sheaf ${\mathcal E}_{\ell}$ on
$X$ is therefore an extension
\[
0 \rightarrow {\mathcal I}_L \rightarrow {\mathcal E}_{\ell} \rightarrow
{\mathcal I}_L(-\ell+1) \rightarrow 0
\]
of rank $1$ ACM sheaves on $X$.

Thus every ACM sheaf on $X$ is layered, and condition (B) is satisfied.

Now the result follows from $(4.3)$.

\bigskip
\noindent
{\bf Theorem 6.2.}  {\em Let $X$ be the singular quadric three-fold in ${\mathbb
P}^4$ that is the cone over a non-singular quadric surface in ${\mathbb P}^3$. 
Two locally Cohen--Macaulay curves on $X$ are in the same Gorenstein biliaison
class if and only if their Rao modules are isomorphic.  In particular, all {\em
ACM} curves are in the same biliaison class.}

\bigskip
\noindent
{\em Proof.} We will show that conditions (B), (C) of $(4.3)$ and (B$''$) of
$(4.7)$ hold, and conclude by applying those two theorems.

Again we use the result of \cite[4.1]{BGS} mentioned in the previous proof, but
in this case we will have to use Kn\"orrer's double branched covers theorem
\cite[2.5]{K}:  if $S_1 = k[x,y,t]$ and $R_1 = S_1/(x^2+y^2)$, then every
indecomposable MCM module on $R_1$ is a direct summand of a module with  matrix
factorization
\[
\left( \begin{pmatrix}
y & \psi \\
\varphi & -y
\end{pmatrix}, \begin{pmatrix}
y & \psi \\
\varphi & -y
\end{pmatrix} \right)
\]
for each matrix factorization $(\varphi,\psi)$ of an indecomposable MCM module
over $R$.

From (i) we obtain the ring $R_1$.  From (ii) we obtain a matrix factorization
\[
\left( \begin{pmatrix}
y & x \\
x & -y
\end{pmatrix}, \begin{pmatrix}
y & x \\
x & -y
\end{pmatrix} \right)\ .
\]
Setting $a = x + iy$, $b = x - iy$, we see this matrix factorization is equivalent
to
\[
\left( \begin{pmatrix}
a & 0 \\
0 & b 
\end{pmatrix}, \begin{pmatrix}
b & 0 \\
0 & a
\end{pmatrix} \right)\ .
\]
The corresponding module is a direct sum of $S_1/(a)$ and $S_1/(b)$.  

From type (iii) we obtain a matrix factorization
\[
\left( \begin{pmatrix}
y & 0 & x & t^{\ell} \\
0 & y & 0 & -x \\
x & t^{\ell} & -y & 0 \\
0 & -x & 0 & -y
\end{pmatrix}, \mbox{ ditto}\right)\ .
\]
Again setting $x^2+y^2 = ab$, we find this matrix factorization is equivalent to
one that corresponds to a direct sum of two modules $M_{\ell}$ and $M'_{\ell}$
with matrix factorizations
\[
\left( \begin{pmatrix}
a & -t^{\ell} \\
0 & b 
\end{pmatrix}, \begin{pmatrix}
b & t^{\ell} \\
0 & a
\end{pmatrix} \right)
\]
and
\[
\left( \begin{pmatrix}
b & -t^{\ell} \\
0 & a
\end{pmatrix}, \begin{pmatrix}
a & t^{\ell} \\
0 & b
\end{pmatrix} \right)\ .
\]

Thus we have shown that indecomposable MCM modules over $R_1$ correspond to
matrix factorizations of $5$ types

\begin{itemize}
\item[(a)] $((x^2+y^2),(1))$
\item[(b)] $((a),(b))$
\item[(c)] $((b),(a))$
\item[(d)] $\left( \begin{pmatrix}
a & -t^{\ell} \\
0 & b
\end{pmatrix}, \begin{pmatrix}
b & t^{\ell} \\
0 & a
\end{pmatrix}\right)$, $\ell = 1,2,\dots$
\item[(e)] $\left( \begin{pmatrix}
b & -t^{\ell} \\
0 & a
\end{pmatrix}, \begin{pmatrix}
a & t^{\ell} \\
0 & b
\end{pmatrix} \right)$, $\ell = 1,2,\dots$.
\end{itemize}

Now we use Kn\"orrer's periodicity theorem \cite[3.1]{K} to pass to the rings
$S_3 = k[x,y,u,v,t]$ and $R_3 = S_3/(x^2+y^2+uv)$.  Preserving the notation $ab =
x^2+y^2$, as in the proof of $(6.1)$ to each of these matrix factorizations
$(\varphi,\psi)$ we obtain matrix factorizations
\[
\left( \begin{pmatrix}
u & \psi \\
\varphi & -v
\end{pmatrix}, \begin{pmatrix}
v & \psi\\
\varphi & -u
\end{pmatrix} \right)
\]
over the ring $S_3$.

From case (a) we obtain the ring $R_3$.  From cases (b) and (c) we obtain modules
isomorphic to the ideals $(a,u)$ and $(a,v)$ in $R_3$.  These correspond to the
two families of planes $D,E$ passing through the vertex of the cone $X$.

From case (d) we obtain a matrix factorization
\[
\left( \begin{pmatrix}
u & 0 & a & -t^{\ell} \\
0 & u & 0 & b \\
b & t^{\ell} & -v & 0 \\
0 & a & 0 & -v
\end{pmatrix}, \begin{pmatrix}
v & 0 & a & -t^{\ell} \\
0 & v & 0 & b \\
b & t^{\ell} & -u & 0 \\
0 & a & 0 & -u
\end{pmatrix} \right)\ .
\]
As in the proof of $(6.1)$, we rearrange rows and columns, and then observe that
the corresponding module $M_{\ell}$ is an extension
\[
0 \rightarrow I_D \rightarrow M_{\ell} \rightarrow I_E(-\ell + 1) \rightarrow 0.
\]
Similarly in case (e) we obtain a module $M'_{\ell}$ that is an extension
\[
0 \rightarrow I_E \rightarrow M'_{\ell} \rightarrow I_D(-\ell+1) \rightarrow 0.
\]

Translating this in terms of the ACM sheaves on $X$, we find that the
indecomposable ACM sheaves on $X$, up to twist, are ${\mathcal O}_X$, ${\mathcal
I}_D$, ${\mathcal I}_E$, and two infinite sequences ${\mathcal
E}_{\ell},{\mathcal E}'_{\ell}$ for $\ell = 1,2,\dots$ that are extensions
\[
0 \rightarrow {\mathcal I}_D \rightarrow {\mathcal E}_{\ell} \rightarrow
{\mathcal I}_E(-\ell+1) \rightarrow 0
\]
\[
0 \rightarrow {\mathcal I}_E \rightarrow {\mathcal E}'_{\ell} \rightarrow
{\mathcal I}_D(-\ell+1) \rightarrow 0.
\]
Thus we see that every ACM sheaf on $X$ is layered, so $X$ satisfies (B) and
(B$''$).  We know that $\mbox{APic } X = {\mathbb Z}^2$, generated by $D,E$ and
that $D+E=H.X$ is arithmetically Gorenstein, and there is an exact sequence
\[
0 \rightarrow {\mathcal I}_D \rightarrow {\mathcal O}_X^2 \rightarrow {\mathcal
I}_E(1) \rightarrow 0
\]
so condition (C) is satisfied.

Now the result follows from $(4.3)$ and $(4.7)$.

\bigskip
\noindent
{\bf Remark 6.3.}  Here is another way to construct the sheaves ${\mathcal
E}_{\ell}$ and ${\mathcal E}'_{\ell}$.  Let $C$ be a plane curve of degree $\ell$
in the plane $D$.  Then $C$ is arithmetically Gorenstein with $\omega_C \cong
{\mathcal O}_C(\ell-3)$.  So by the Serre construction there is an extension
\[
0 \rightarrow {\mathcal O}_X(-\ell) \rightarrow {\mathcal E} \rightarrow
{\mathcal I}_C \rightarrow 0,
\]
and ${\mathcal E}$ will be a rank $2$ ACM sheaf.  Since $C \subseteq D$, there is
an inclusion $0 \rightarrow {\mathcal I}_D \rightarrow {\mathcal E}$, and chasing
the corresponding diagram one finds the cokernel is ${\mathcal I}_E(-\ell+1)$. 
Thus ${\mathcal E}$ is the sheaf ${\mathcal E}_{\ell}$ described above.

\bigskip
\noindent
{\bf Example 6.4.}  As an application of $(6.2)$ we give a new proof and
strengthening of a theorem of Lesperance.  Let $D_1,D_2$ be two planes in the
family of planes $D$, meeting at the singular point $P = (0,0,0,0,1)$ of $X$.  Let
$C_1
\subseteq D_1$ and $C_2 \subseteq D_2$ be curves of degrees $d,t$ respectively,
not containing the point $P$, with $2 \le d \le t$.  Lesperance calls the union
$C = C_1 \cup C_2$ a curve of type $(P,d,t)$, and shows that its Rao module is
$M_d = S_3/(x,y,u,v,t^d)$.  Then he shows that two curves $C$
and $C'$ of types $(P,d,t)$ and $(P,d,s)$ with $2 \le d \le t,s$ are in the same
$G$-liaison class \cite[4.9]{L1}, \cite[3.6]{L2}.  Since they have the same Rao
module, we obtain a new proof of this result from $(6.2)$, which shows they are
also in the same $G$-biliaison class.

Next, Lesperance considers a curve $C' = C_1 \cup C'_2$ as above, where $C_1
\subseteq D_1$ does not contain $P$, but where $C'_2 \subseteq D_2$ does contain
$P$ as a point of multiplicity $e$.  This time we assume $d = \mbox{deg } D_1 \ge
2$, but $t \ge e \ge 1$ only.  This he calls a curve of type $(P,d,t,e)$, and
shows it has the same Rao module $M_d$.  He shows that if $t \ge d + e$, then
curves of type $(P,d,t)$ and $(P,d,t,e)$ are in the same even $G$-liaison class
\cite[4.15]{L1}, but he leaves open the question if $d = t$, and asks, for
example if curves of types $(P,2,2)$ and $(P,2,2,1)$ are $G$-linked.  Again,
since they have the same Rao modules, our theorem $(6.2)$ shows that any curves
of types $(P,d,t)$ and $(P,d,t,e)$ are in the same $G$-biliaison class.  This
answers the question \cite[4.14]{L1}.

Note that our $G$-biliaisons take place inside a fixed singular quadric
hypersurface $X$.  But since any three skew lines in ${\mathbb P}^3$ are
contained in a non-singular quadric surface, it is easy to see that if
$D'_1,D'_2$ are any other two planes in ${\mathbb P}^4$ meeting at the same point
$P$, then we can make $G$-biliaisons on different singular quadric hypersurfaces
to relate curves of types $(P,d,t)$ or $(P,d,t,e)$ in $D_1$ and $D_2$ to those in
$D'_1$ and $D'_2$.

\bigskip
\noindent
{\bf Example 6.5.}  Lesperance also gives an example \cite[5.7]{L1},
\cite[4.6]{L2} of minimal curves in the same even $G$-liaison class, having the
same degree and genus, that do not belong to the same irreducible family.  His
example consists of a curve $C$ of type $(P,2,2,1)$ on the one hand, and the
disjoint union $D$ of a line $L$ and a twisted cubic curve $Y$, where $L$ meets
the ${\mathbb P}^3$ containing $Y$ at the point $P$ on the other hand.  Both have
degree
$4$,
$p_a = -1$, and Rao module $M_2$.

We would like to point out that we can put curves of type $C$ and $D$ on a
singular quadric $3$-fold $X$, so that they are in the same $G$-biliaison class
on $X$.  This will show that minimal curves in a $G$-biliaison class on $X$, of
the same degree and genus, need not form an irreducible family.

To do this, let $Q$ be a non-singular quadric surface in ${\mathbb P}^3$.  Take a
point $O \in {\mathbb P}^4{\backslash}{\mathbb P}^3$, and let $X$ be the cone
over $Q$ with vertex $O$.  Now take $Y \subseteq Q$ a twisted cubic curve, take
$P \in Q{\backslash}Y$, and let $L$ be the line joining $O$ to $P$.  Then $D = L
\cup Y$ has Rao module $M_2$ based at $P$.  Take a plane $\Lambda_1$ in ${\mathbb
P}^3$ cutting $Q$ in a conic $C_1$ containing $P$.  Let $M$ be a line in $Q$,
passing through $P$, not in the plane $\Lambda_1$, and let $\Lambda_2$ be the
plane spanned by $M$ and $O$.  Let $C_2 \subseteq \Lambda_2$ be a conic not
passing through $P$.  Then $C = C_1 \cup C_2$ is a curve of type $(P,2,2,1)$. 
Since $C$ and $D$ have the same Rao module, they are in the same $G$-biliaison
class on $X$, by $(6.2)$.

\end{document}